\documentclass[a4paper,11pt]{article}
\usepackage{xypic}
\usepackage{tikz}
\usepackage{amsmath,amssymb,amsthm,bm,mathtools}
\usepackage{mathrsfs}
\usepackage{makeidx}
\makeindex

\def\lim{\mathop{\rm lim}}

\def\pf{{\it Proof:}~}

\newtheorem{theorem}{Theorem}
\newtheorem{lemma}[theorem]{Lemma}
\newtheorem{proposition}[theorem]{Proposition}
\newtheorem{definition}[theorem]{Definition}

\newtheorem{remark}[theorem]{Remark}

\newcommand{\ud}{\, \mathrm{d}}

\title{\bf On the limiting velocity of random walks in mixing random environment}
\author{Xiaoqin Guo\thanks{School of Mathematics, University of Minnesota,
206 Church St SE, Minneapolis, MN 55455. Research supported by NSF
grant DMS-0804133.}}
\date{June 28, 2011. Revised September 9, 2012}
\begin{document}
\maketitle
\begin{abstract}
We consider random walks in strong-mixing random Gibbsian environments in $\mathbb{Z}^d, d\ge 2$.
Based on regeneration arguments, we will first provide an alternative proof of Rassoul-Agha's conditional law of large numbers (CLLN) for mixing environment \cite{R-A3}.
Then, using coupling techniques, we show that there is at most one nonzero limiting velocity in high dimensions ($d\ge 5$).
\end{abstract}

\section{Introduction}
An \textit{environment} is an element $\omega=\{\omega(x,e)\}_{x\in\mathbb{Z}^d, |e|=1}$ 
of $\Omega=\mathcal{M}^{\mathbb{Z}^d}$,
where $\mathcal{M}$ is the space of probability measures on $\{e\in \mathbb{Z}^d:|e|=1\}$ and
$|\cdot|$ denotes the Euclidean norm.
The random walk in the
environment $\omega\in\Omega$
started at $x$ is the canonical Markov
chain $(X_n)$ on $(\mathbb{Z}^d)^\mathbb{N}$,
with state space $\mathbb{Z}^d$ and law $P_\omega^x$ specified by
\begin{eqnarray*}
&&P_\omega^x\{X_0=x\}=1,\\
&&P_\omega^x\{X_{n+1}=y+e | X_n=y\}=\omega(y, e), \quad e\in\mathbb{Z}^d, |e|=1.
\end{eqnarray*}

Let $P$ be a stationary (with respect to the shifts in $\mathbb{Z}^d$) probability measure on $\Omega$. The joint law of the environment and the walks is denoted by $\textbf{P}^x=P\otimes P_\omega^x$. We also write $\textbf{P}^o$ as $\textbf{P}$, where
$o$ denotes the origin. We say that the random environment is \textit{iid} if $P$ is a product measure. We say that $P$ is \textit{uniformly elliptic} if there is a constant $\kappa\in (0,1/2d)$ such that
$P$-almost surely,
\[
\omega(o,e)>\kappa \text{ for all $e\in\mathbb{Z}^d$ with $|e|=1$}.
\]

For any vector $\ell\in S^{d-1}$, we let 
\[
A_{\ell}=\{\lim_{n\to\infty}X_n\cdot\ell=\infty\}.
\]

In recent years, much progress has been made in the study of the limiting velocity
$\lim_{n\to\infty}X_n/n$ of random
walks in random iid environment, see \cite{ZO} for a survey. 
For one-dimensional RWRE, the law of large numbers (LLN) is well known (see \cite{So}).
For $d\ge 2$, a conditional law of large numbers (CLLN) is proved in \cite{SZ, Ze} (see \cite[Theorem 3.2.2]{ZO} for the full version), which states that $\mathbf{P}$-almost surely, for any direction $\ell$,
\[
\lim_{n\to\infty}\frac{X_n\cdot\ell}{n}=v_\ell 1_{A_\ell}-v_{-\ell}1_{A_{-\ell}} \tag{CLLN}
\]
for some deterministic constants $v_\ell$ and $v_{-\ell}$ (we set $v_\ell=0$ if $\mathbf{P}(A_\ell)=0$). 
Moreover, for $d=2$, the LLN follows from combining the CLLN and Zerner and Merkl's $0$-$1$ law \cite{ZM} for two-dimensional RWRE: for any direction $\ell$,
\[
\mathbf{P}(A_\ell)\in\{0,1\}.
\]
When $d\ge 3$, the $0$-$1$ law and the LLN are among the main open questions in the study of RWRE. Nevertheless, in high dimension ($d\ge 5$), Berger \cite{Be} showed
that the limiting velocity can take at most one non-zero value, i.e.,
\begin{equation}\label{berger}
v_\ell v_{-\ell}=0.
\end{equation}

The purpose of this paper is to extend the CLLN and Berger's result
\eqref{berger} to the case when the environments on different sites are allowed
to be dependent. Of special interest is the environment that is produced by a
Gibbsian particle system (which we call the \textit{Gibbsian environment}) and
satisfies Dobrushin-Shlosman's strong-mixing condition IIIc in \cite[page
378]{DS}, see \cite{R-A1,R-A2,CZ1,CZ2, R-A3} for related works.  For the
definition of the Gibbsian environment and the strong-mixing condition
\cite[(6.1)]{R-A1}, see \cite[pages 1454-1455]{R-A1}. An important feature of
this model is that the influence of the environments in remote locations decays
exponentially as the distance grows.

In \cite{R-A1}, assuming a ballisticity condition (Kalikow's condition) which implies that the 
event of escape in a direction has probability $1$,  Rassoul-Agha proved the LLN for the strong-mixing Gibbsian environment, using the invariant measure of the ``environment viewed from the point of view of the particle" process. 
In \cite{R-A3}, Rassoul-Agha also obtained a CLLN for the strong-mixing Gibbsian environment, under an analyticity condition (see Hypothesis (M) in \cite{R-A3}).
Comets and Zeitouni proved the LLN for environments with a weaker cone-mixing assumption ($\mathcal{A}1$) in \cite{CZ1}, but under some conditions about ballisticity and the uniform integrability of the regeneration times (see ($\mathcal{A}5$) in \cite{CZ1}).

Our first purpose is to prove the CLLN for random walks in the strong-mixing
Gibbsian environment. Display \eqref{LVclln} in Theorem~\ref{thm2} is a minor
extension of Rassoul-Agha's CLLN
in \cite{R-A3}, in which he assumes slightly more than strong-mixing.
Yet, our proof is very different from the proof in 
\cite{R-A3} , which is based on a large deviation principle in \cite{R-A2}. The
main contribution of our proof of \eqref{LVclln} is a new definition of the
regeneration structure, which enables us to divide a random path in the mixing
environment into ``almost iid" parts.
With this regeneration structure, we will use the ``$\epsilon$-coins" introduced
in \cite{CZ1} and coupling arguments to prove the CLLN. This regeneration
structure will also be used in the proof of \eqref{LVunique}.

\begin{figure}
\centering
\includegraphics[width=0.7\textwidth]{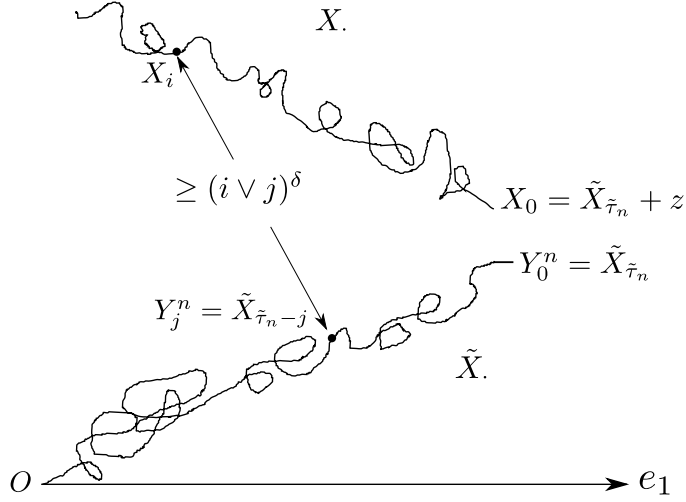}
\caption{$X_\cdot\in A(Y_\cdot^n, z)$. 
When $i\vee j>|z|/3$, the distance between $Y_j^n$ of the ``backward path" and $X_i$ is at least $(i\vee j)^\delta$.}
\label{fig:1}
\end{figure}

Our second main result \eqref{LVunique} is an extension of Berger's result
\eqref{berger} from the iid case to the strong-mixing case. 
In \cite{Be}, assuming that $\mathbf{P}(A_\ell)>0$ for a direction $\ell$, Berger coupled the iid environment $\omega$ with a transient (in the direction $\ell$) environment $\tilde\omega$ and a ``backward path", such that $\tilde\omega$ and $\omega$ coincide in the locations off the path. 
Using heat kernel estimates for random walks with iid increments, he showed that if $v_\ell v_{-\ell}>0$ and $d\ge 5$, then with positive probability, the random walks in $\tilde\omega$ is transient to the $-\ell$ direction without intersecting the backward path, which contradicts $\tilde{\omega}$ being transient in the direction $\ell$.
The difficulties in applying this argument to mixing environments are that the regeneration slabs are not
iid, and that unlike the iid case, the environments visited by two disjoint paths are not independent. To overcome these difficulties,
we will construct an environment (along with a path) that is ``very transient" in $\ell$,  and show that the ballistic walks in the opposite direction ($-\ell$) will move further and further away from the given path (see Figure \ref{fig:1} in Section \ref{secunique}). The key ingredient here is a heat kernel estimate, which we will obtain in Section \ref{sechke} using coupling arguments. 

We now describe our main results. Recall first the definition of an $r$-Markov environment (see \cite{CZ2}). 
\begin{definition}\label{def1}
For $r\ge 1$, let $\partial_r V=\{x\in\mathbb{Z}^d\setminus V: d(x, V)\le r\}$ be the $r$-boundary of $V\subset\mathbb{Z}^d$.
A random environment $(P,\Omega)$ on $\mathbb{Z}^d$ is called $r$-Markov if 
for any finite $V\subset\mathbb{Z}^d$,
\[
P\big((\omega_x)_{x\in V}\in \cdot|\mathcal{F}_{V^c}\big)
=P\big((\omega_x)_{x\in V}\in \cdot|\mathcal{F}_{\partial_r V}\big), \text{ $P$-a.s.,}
\]
where $d(\cdot,\cdot)$ denotes the $l^1$-distance and $\mathcal{F}_{\Lambda}:=\sigma(\omega_x:x\in\Lambda)$.
\end{definition}
We say that an $r$-Markov environment $P$ satisfies condition ($*$) if there exist 
constants $\gamma , C<\infty$ such that for all finite subsets $\Delta\subset V\subset\mathbb{Z}^d$ with $d(\Delta,V^c)\ge r$, and $A\subset V^c$,
\[
\frac{\ud P\big((\omega_x)_{x\in\Delta}\in\cdot|\eta\big)}
{ \ud P\big((\omega_x)_{x\in\Delta}\in\cdot|\eta'\big)}
\le 
\exp{(C\sum_{x\in A,y\in\Delta}e^{-\gamma d(x,y)})}\tag{$*$}
\]
for $P$-almost all pairs of configurations $\eta,\eta'\in\mathcal{M}^{V^c}$ which agree on $V^c\setminus A$.
Here 
\[
P\big((\omega_x)_{x\in\Delta}\in\cdot|\eta\big)
:=P\big((\omega_x)_{x\in\Delta}\in \cdot|\mathcal{F}_{V^c}\big)\big|_{(\omega_x)_{x\in V^c}=\eta}.
\]
We remark that $r$ and $\gamma$ are used as parameters of the environment throughout the article.

By Lemma 9 in \cite{R-A1}, the strong-mixing Gibbsian environment satisfies ($*$).
Obviously, every finite-range dependent environment also satisfies ($*$).

Our main theorem is:
\begin{theorem}\label{thm2}
Assume that $P$ is uniformly elliptic and satisfies \emph{($*$)}.
Then there exist two deterministic constants $v_+, v_-\ge 0$
and a vector $\ell$ such that
\begin{equation}\label{LVclln}
\lim_{n\to\infty} \frac{X_n}{n}=v_+\ell 1_{A_\ell}-v_-\ell 1_{A_{-\ell}},
\end{equation}
and $v_+=v_-=0$ if $\mathbf{P}(A_\ell\cup A_{-\ell})<1$.
Moreover, if $d\ge 5$, then
there is at most one non-zero velocity. That is, 
\begin{equation}\label{LVunique}
v_+ v_-=0.
\end{equation}
\end{theorem}
We remark here that for the finite-range dependent
case, the CLLN is proved in \cite{ZO}.

The structure of this paper is as follows. In Section \ref{seccomb}, we prove a refined 
version of \cite[Lemma 3]{Ze}. With this combinatorial result, we will prove the CLLN in Section \ref{seclln}, using coupling arguments. In Section \ref{sechke}, using coupling, we obtain
heat kernel estimates, which is later used in
Section \ref{secunique} to show the uniqueness of the non-zero limiting velocity. 

Throughout the paper, we assume that the environment is uniformly elliptic and satisfies $(*)$. We use $c, C$ to denote finite positive constants that depend only on 
the dimension $d$ and the environment measure $P$ (and implicitly, on the parameters $\kappa,r$ and $\gamma$ of the environment). They may differ from line to line. 
We denote by $c_1,c_2,\ldots$ positive constants which are fixed throughout, and which depend only on $d$ and the measure $P$. Let $\{e_1,\ldots,e_d\}$ be the natural basis of $\mathbb{Z}^d$.

\section{A combinatorial lemma and its consequences}\label{seccomb}
In this section we consider the case that $\mathbf{P}(\varlimsup_{n\to\infty}X_n\cdot e_1/n>0)>0$. We will adapt the arguments in \cite{Ze} and prove that with positive probability,
the number of visits to the $i$-th level $\mathcal{H}_i=\mathcal{H}_i(X_0):=\{x:x\cdot e_1=X_0\cdot e_1+ i\}$ grows
slower than $Ci^2$.
An important ingredient of the proof is a refinement of a combinatorial lemma of Zerner \cite[Lemma 3]{Ze} about deterministic paths.

We say that a sequence $\{x_i\}_{i=0}^{k-1}\in (\mathbb{Z}^d)^{k}$, $2\le k\le\infty$, is a \textit{path} if 
$|x_i-x_{i-1}|=1$ for $i=1,\cdots, k-1$. For $i\ge 0$ and an infinite path $X_\cdot=\{X_n\}_{n=0}^\infty$ such that $\sup_n X_n\cdot e_1=\infty$, let
\[T_i=\inf\{n\ge 0: X_n\in\mathcal{H}_i\}.\]
For $0\le i<j$ and $k\ge 1$, let $T_{i,j}^1:=T_i$ and define
recursively
\[
T_{i,j}^{k+1}=\inf\{n\ge T_{i,j}^k: X_n\in\mathcal{H}_i \text{ and } n<T_j\}\in \mathbb{N}\cup  \{\infty\}.
\]
That is, $T_{i,j}^k$ is the time of the $k$-th visit to $\mathcal{H}_i$ before hitting
$\mathcal{H}_j$. Let
\[
N_{i,j}=\sup\{k: T_{i,j}^k<\infty\}
\]
be the total number of visits to $\mathcal{H}_i$ before hitting
$\mathcal{H}_j$.

As in \cite{Ze}, for $i\ge 0, l\ge 1$, let
\[
h_{i,l}=T_{i,i+l}^{N_{i,i+l}}-T_i
\]
denote the time spent between the first and the last visits to $\mathcal{H}_i$ before hitting $\mathcal{H}_{i+l}$. 
For $m,M, a\ge 0$ and $l\ge 1$, set
\[
H_{m,l}=\sum_{i=0}^{l-1}N_{m+i,m+l}/(i+1)^2
\]
and
\[
E_{M,l}(a)=\frac{\#\{0\le m\le M: h_{m,l}\le a \text{ and } H_{m,l}\le a\}}{M+1}.
\]
Note that $E_{M,l}(a)$ decreases in $l$ and increases in $a$.

The following lemma is
a minor adaptation of \cite[Lemma 3]{Ze}.
\begin{lemma}\label{l5}
For any path $X_\cdot$ with $\varlimsup_{n\to\infty}X_n\cdot e_1/n>0$, 
\begin{equation}\label{e27}
\sup_{a\ge 0}\inf_{l\ge 1}\varlimsup_{M\to\infty}E_{M,l}(a)>0.
\end{equation}
\end{lemma}
\pf
Since $\varlimsup_{n\to\infty}n/T_n=\varlimsup_{n\to\infty}X_n\cdot e_1/n>0$,
there exist an increasing sequence $(n_k)_{k=0}^\infty$ and $\delta<\infty$ such that
\[
T_{n_k}<\delta n_k \text{ for all }k.
\]
Thus for any $m$ such that $n_k/2\le m\le n_k$,
\begin{equation}\label{*17}
T_m\le 2\delta m.
\end{equation}
Set $M_k=\lceil n_k/2\rceil$. Then for all $k$ and $1<l<\lfloor n_k/2 \rfloor$,
\begin{align}\label{e28}
\sum_{m=0}^{M_k} H_{m,l}
&=\sum_{i=0}^{l-1}\Big(\sum_{m=0}^{M_k}N_{m+i,m+l}\Big)/(i+1)^2\nonumber\\
&\le \sum_{i=0}^{l-1} T_{M_k+l}/(i+1)^2
\stackrel{(\ref{*17})}{\le} 4\delta (M_k+l).
\end{align}

By the same argument as in Page 193-194 of \cite{Ze}, we will show that there exist
constants $c_1, c_2>0$ such that
\begin{equation}\label{*18}
\inf_{l\ge 1}\varlimsup_{k\to\infty}
\frac{\#\{0\le m\le M_k: h_{m,l}\le c_1 \}}{M_k+1}>c_2.
\end{equation}
Indeed, if (\ref{*18}) fails, 
then for any $u>0$, 
\[
	\varlimsup_{k\to\infty}\dfrac{\#\{0\le m\le M_k, h_{m,l}\le u\}}{M_k+1}\longrightarrow 0 
\]
as $l\to\infty$ (note that the right side is decreasing in $l$). Hence,
one can find a sequence $(l_i)_{i\ge 0}$ with $l_{i+1}>l_i, l_0=0,$
such that for all $i\ge 0$,
\begin{equation}\label{*19}
\varlimsup_{k\to\infty}\dfrac{\#\{0\le m\le M_k, h_{m,l_{i+1}}\le 6\delta l_i\}}{M_k+1}<\frac{1}{3}.
\end{equation}
On the other hand, for $i\ge 0$
\begin{align}
&\varlimsup_{k\to\infty}\dfrac{\#\{0\le m\le M_k, h_{m,l_i}\ge 6\delta l_i\}}{M_k+1}\nonumber\\
&\le 
\varlimsup_{k\to\infty}\frac{1}{(M_k+1)6\delta l_i}\sum_{m=0}^{M_k}(T_{m+l_i}-T_m)\nonumber\\
&\le 
\varlimsup_{k\to\infty}\frac{l_i T_{M_k+l_i}}{6\delta l_i(M_k+1)}
\stackrel{(\ref{*17})}{\le}\frac{1}{3}. \label{*20}
\end{align}
By (\ref{*19}) and (\ref{*20}) 
, for any $i\ge 0$,
\begin{equation}\label{*21}
\varlimsup_{k\to\infty}\dfrac{\#\{0\le m\le M_k, h_{m,l_{i+1}}> h_{m,l_i}\}}{M_k+1}
\ge \frac{1}{3}.
\end{equation}
Therefore, for any $j\ge 1$, noting that 
\[
\sum_{i=0}^{j-1}1_{h_{m,l_{i+1}}>h_{m,l_i}}\le N_{m,m+l_j}\le H_{m,l_j},
\]
 we have 
\begin{align*}
\frac{j}{3}&\stackrel{(\ref{*21})}{\le}
\varlimsup_{k\to\infty}
\sum_{i=0}^{j-1}\dfrac{\#\{0\le m\le M_k, h_{m,l_{i+1}}> h_{m,l_i}\}}{M_k+1}\\
&\le\varlimsup_{k\to\infty}
\frac{1}{M_k+1}\sum_{m=0}^{M_k}H_{m,l_j}\stackrel{(\ref{e28})}{\le} 4\delta,
\end{align*}
which is a contradiction if $j$ is large. This proves (\ref{*18}).

It follows from (\ref{*18}) that, for any $l\ge 1$, there is a subsequence
$(M'_k)$ of $(M_k)$ such that
\[
\frac{\#\{0\le m\le M'_k: h_{m,l}\le c_1 \}}{M'_k+1}>c_2
\]
for all $k$. 
Letting $c_3=9\delta/c_2$, we have that when $k$ is large enough,
\[
\frac{1}{M'_k+1}\sum_{m=0}^{M'_k}1_{h_{m,l}\le c_1, H_{m,l}>c_3}
\le 
\frac{1}{c_3(M'_k+1)}\sum_{m=0}^{M'_k}H_{m,l}
\stackrel{(\ref{e28})}{\le}
\frac{c_2}{2}.
\]
Hence for any $l>1$ and large $k$,
\begin{align*}
E_{M_k',l}(c_1\vee c_3)
&\ge \frac{1}{M'_k+1}\sum_{m=0}^{M'_k}1_{h_{m,l}\le c_1,H_{m,l}\le c_3}\\
&= 
\frac{1}{M'_k+1}\sum_{m=0}^{M'_k}(1_{h_{m,l}\le c_1}-1_{h_{m,l}\le c_1,H_{m,l}> c_3})
\ge \frac{c_2}{2}.
\end{align*}
This shows the lemma, and what is more, with explicit constants.\qed\\

For $i\ge 0$, let $N_i=\lim_{j\to\infty} N_{i,j}$ denote the total number of visits to $\mathcal{H}_i$.
With Lemma \ref{l5}, one can deduce that with positive probability, $N_i\le C(i+1)^2$ for all $i\ge 0$:
\begin{theorem}\label{thm3}
If $\mathbf{P}(\varlimsup_{n\to\infty}X_n\cdot e_1/n>0)>0$, then there exists a constant
$c_5$ such that
\[\mathbf{P}(R=\infty)>0,\]
where $R$ is the stopping time defined by
\begin{align*}
R&=R_{e_1}(X_\cdot, c_5)\\
&:=
\inf\{n\ge 0: \sum_{i=0}^n 1_{X_i\in\mathcal{H}_j}>c_5(j+1)^2 \text{ for some }j\ge 0\}\wedge D,
\end{align*}
and $D:=\inf\{n\ge 1: X_n\cdot e_1\le X_0\cdot e_1\}$.
\end{theorem}

\begin{figure}[h]
\centering
\includegraphics[width=0.7\textwidth]{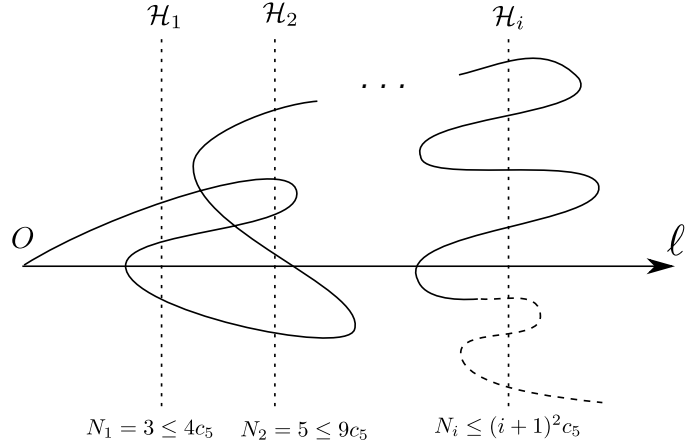}
\caption{On $\{R=\infty\}$, the path visits the $i$-th level no more than
$c_5(i+1)^2$ times.}
\label{intersection}
\end{figure}

Note that for any $L>0$ and a path $(X_i)_{i=0}^\infty$ with $X_0=o$,
\begin{align}\label{e32}
\sum_{\substack{y: y\cdot e_1\le -L\\0\le i\le R}}e^{-\gamma d(y,X_i)}
&\le
\sum_{j=0}^\infty (\#\text{visits to $\mathcal{H}_j$ before time $R$})e^{-\gamma (j+L)}\nonumber\\
&\le 
C\sum_{j=0}^\infty c_5(j+1)^2e^{-\gamma(j+L)}\le Ce^{-\gamma L}.
\end{align}
Hence on the event $\{R=\infty\}$, by (\ref{e32}) and $(*)$,
 the trajectory $(X_i)_{i=0}^\infty$ is ``almost 
independent" with the environments $\{\omega_x:x\cdot e_1\le -L\}$ when $L$ is
large. See Figure~\ref{intersection}. 
This fact will be used in our definition of the regeneration times in the Section \ref{seclln}.

To prove Theorem \ref{thm3}, we need the following lemma. 
Recall that $r,\gamma$ are parameters of the environment measure $P$.
Let $S$ be a countable set of finite paths.
With abuse of notation, we
also use $S$ as the synonym for the event
\begin{equation}\label{2e*}
\bigcup_{(x_i)_{i=0}^N\in S}\{X_i=x_i \text{ for }0\le i\le N\}.
\end{equation}
\begin{lemma}\label{c2}
Let $a>0$ and $A\subset\Lambda\subset\mathbb{Z}^d$.
Suppose $S\neq\emptyset$ is a countable set of finite paths 
$x_\cdot=(x_i)_{i=0}^N, N<\infty$ that satisfy $d(x_\cdot, \Lambda)\ge r$ and
\[
\sum_{y\in A, 0\le i\le N}e^{-\gamma d(y,x_i)}\le a.
\]
Then, $P$-almost surely,
\begin{equation}\label{e31}
	\exp(-Ca)\le\frac{E_P [P_\omega(S)|\omega_x: x\in\Lambda]}{E_P [P_\omega(S)|\omega_x: x\in\Lambda\setminus A]}\le\exp(Ca).
\end{equation}
\end{lemma}

\pf
We shall first show that for any $(x_i)_{i=0}^N\in S$, $P$-almost surely,
\begin{align}\label{e40}&E_P[P_\omega(X_i=x_i,0\le i\le N)|\omega_y:y\in\Lambda]\nonumber\\&\le \exp(Ca)E_P[P_\omega(X_i=x_i,0\le i\le N)|\omega_y:y\in\Lambda\setminus A].\end{align}
Note that when $\Lambda^c$ is a finite subset of $\mathbb{Z}^d$,  (\ref{e40})
is an easy consequence of $(*)$. For general $\Lambda$, we let
\[
\Lambda_n=\Lambda\cup\{x:|x|\ge n\}.
\]
When $n$ is sufficiently big, $(*)$ implies that 
\begin{equation*}
	\frac{E_P[P_\omega(X_i=x_i,0\le i\le N)|\omega_y:y\in\Lambda_n]}{E_P[P_\omega(X_i=x_i,0\le i\le N)|\omega_y:y\in\Lambda_n\setminus A]}\le\exp(Ca).
\end{equation*}
Since $\Lambda_n\downarrow \Lambda$ as $n\to\infty$, 
(\ref{e40}) follows by taking $n\to\infty$ in the above inequality.

Summing over all $(x_i)_{i=0}^N\in S$ on both sides of (\ref{e40}), we conclude that
$P$-almost surely,
\[
E_P[P_\omega(S)|\omega_y:y\in\Lambda]
\le \exp(Ca)E_P[P_\omega(S)|\omega_y:y\in\Lambda\setminus A]. 
\]
The upper bound of (\ref{e31}) is proved. The lower bound follows likewise.\qed

Now we can prove the theorem. Our proof is a modification
of the proof of Theorem 1 in \cite{Ze}:\\

\noindent\textit{Proof of Theorem \ref{thm3}:}
It follows by Lemma \ref{l5} that there exists a constant $c_4>0$ such that
\begin{equation}\label{*22}
\mathbf{P}(\inf_{l\ge 1}\varlimsup_{M\to\infty}E_{M,l}(c_4)>0)>0.
\end{equation}
For $l>r$, $k\ge 0$ and $z\in\mathbb{Z}^d$ with $z\cdot e_1=r$, let $B_{m,l}(z,k,c)$ denote the event
\[
\{N_{m+r,m+l}=k,X_{T_{m+r,m+l}^k}=X_{T_m}+z,H_{m+r,l-r}\le c\}.
\]

Note that on the event $\{h_{m,l}\le c_4\text{ and }H_{m,l}\le c_4\}$, we have
\begin{align*}
T_{m+r,m+l}^{N_{m+r,m+l}}-T_m
&\le h_{m,l}+\sum_{i=0}^r N_{m+i,m+l}\\
&\le c_4+\sum_{i=0}^r (i+1)^2c_4\le (1+r)^3c_4,
\shortintertext{and}
H_{m+r,l-r}
&\le \sum_{i=0}^{l-r-1}(r+1)^2N_{m+r+i,m+l}/(r+i+1)^2\\
&\le (r+1)^2c_4=:c_5.
\end{align*}
Hence
$
\{h_{m,l}\le c_4\text{ and }H_{m,l}\le c_4\}\subset \bigcup_{|z|,k\le (r+1)^3c_4}B_{m,l}(z,k,c_5),
$
and
\[
\lim_{l\to\infty}\varlimsup_{M\to\infty}E_{M,l}(c_4)
\le 
\sum_{|z|,k\le (r+1)^3c_4}
\varlimsup_{l\to\infty}\varlimsup_{M\to\infty}\frac{1}{M+1}
\sum_{m=0}^M 1_{B_{m,l}(z,k,c_5)}.
\]
Thus by (\ref{*22}), for some $k_0$ and $z_0$ with $z_0\cdot e_1=r$,
\begin{equation}\label{e29}
\mathbf{P}(\varlimsup_{l\to\infty}\varlimsup_{M\to\infty}
\frac{1}{M+1}\sum_{m=0}^M 1_{B_{m,l}(z_0,k_0,c_5)}>0)>0.
\end{equation}
In what follows, we write $B_{m,l}(z_0,k_0,c_5)$ simply as
$B_{m,l}$. 

For any $l>r$ and any fixed $i\le l-1$, let $m_j=m_j(l,i):=i+jl$, i.e. $(m_j)_{j\ge 0}$ is
the class of residues of $i(\text{mod }l)$.
Now take any $j\in \mathbb{N}$.	Observe that for any event $E=\{1_{B_{m_{j-1},l}}=\cdot,\ldots,1_{B_{m_0,l}}=\cdot\}$
and $x\in \mathcal{H}_{m_j}$,
\begin{align}\label{*1}
\MoveEqLeft P_\omega(\{X_{T_{m_j}}=x\}\cap E\cap B_{m_j,l})\\
&\le 
P_\omega(\{X_{T_{m_j}}=x\}\cap E)
P_\omega^{x+z_0}(D>T_{l-r},H_{0,l-r}\le c_5).\nonumber
\end{align} 
Moreover, for any $x\in \mathcal{H}_{m_j}$, there exists a countable set $S$ of finite paths $(x_i)_{i=0}^N$ that satisfy
$m_j+r\le x_i\cdot e_1\le m_j+l$ and $\#\{k\le N: x_k\in\mathcal{H}_i(x_0)\}\le c_5(i+1)^2$ for
$0\le i\le N$, such that
\begin{align*}
&\{X_0=x+z_0, D>T_{l-r},H_{0,l-r}\le c_5\}\\
&=\cup_{(x_i)_{i=0}^N\in S}\{X_i=x_i
\text{ for }0\le i\le N\}.
\end{align*}
Noting that (by the same argument as in (\ref{e32})) for any $(x_i)_{i=0}^N\in S$,
\[
\sum_{\substack{y:y\cdot e_1\le m_j\\i\le N}}e^{-\gamma d(y,x_i)}\le Ce^{-\gamma r},
\]
by Lemma \ref{c2} we have
\begin{align*}
&E_P[P_\omega^{x+z_0}(D>T_{l-r},H_{0,l-r}\le c_5)|\omega_y:y\cdot e_1\le m_j]\\
&\le 
\exp{(Ce^{-\gamma r})}\mathbf{P}(D>T_{l-r},H_{0,l-r}\le c_5).
\end{align*}
Thus for $j\ge 0$ and $l>r$,
\begin{align*}
&\mathbf{P}(E\cap B_{m_j,l})\\
&\stackrel{(\ref{*1})}{\le}
\sum_{x\in\mathcal{H}_{m_j}}
E_P \big[P_\omega(\{X_{T_{m_j}}=x\}\cap E)P_\omega^{x+z_0}(D>T_{l-r},H_{0,l-r}\le c_5)\big]\\
&\le 
\exp{(Ce^{-\gamma r})}
\sum_{x\in\mathcal{H}_{m_j}}
\mathbf{P}(\{X_{T_{m_j}}=x\}\cap E)\mathbf{P}(D>T_{l-r},H_{0,l-r}\le c_5)\\
&=
C\mathbf{P}(E)
\mathbf{P}(D>T_{l-r},H_{0,l-r}\le c_5).
\end{align*}
Hence, for any $j\ge 0$ and $l>r$,
\begin{equation*}
\mathbf{P}(1_{B_{m_j,l}}=1|1_{B_{m_{j-1},l}},\ldots,1_{B_{m_0,l}})
\le 
C\mathbf{P}(D>T_{l-r},H_{0,l-r}\le c_5),
\end{equation*}
which implies that $\mathbf{P}$-almost surely,
\begin{equation}\label{e30}
\varlimsup_{n\to\infty}
\frac{1}{n}\sum_{j=0}^{n-1} 1_{B_{m_j,l}}
\le C\mathbf{P}(D>T_{l-r},H_{0,l-r}\le c_5).
\end{equation}
Therefore, $\mathbf{P}$-almost surely,
\begin{align*}
\varlimsup_{l\to\infty}\varlimsup_{M\to\infty}
\frac{1}{M+1}\sum_{m=0}^M 1_{B_{m,l}}
&\le \varlimsup_{l\to\infty}\frac{1}{l}\sum_{i=0}^{l-1}
\varlimsup_{M\to\infty}\frac{l}{M+1}
\sum_{\substack{0\le m\le M\\m\text{ mod }l=i}} 1_{B_{m,l}}\\
&\stackrel{(\ref{e30})}{\le} \lim_{l\to\infty}
C\mathbf{P}(D>T_{l-r},H_{0,l-r}\le c_5)\\
&=C\mathbf{P}(D=\infty, \sum_{i=0}^\infty N_i/(i+1)^2\le c_5).
\end{align*}
This and (\ref{e29}) yield
$\mathbf{P}(D=\infty, \sum_{i=0}^\infty N_i/(i+1)^2\le c_5)>0$. 
The theorem follows.\qed

\section{The conditional law of large numbers}\label{seclln}
In this section we will prove the conditional law of large numbers
\eqref{LVclln}, using
regeneration times and coupling. 
Given the dependence structure of the environment, we want 
to define regeneration times in such a way that what happens after a regeneration time has little 
dependence on the past. To this end, we will use the ``$\epsilon$-coins" trick introduced in \cite{CZ1} and the stopping time $R$ to define the regeneration times. 
Intuitively, at a regeneration time, the past and the future movements have nice properties. That
is, the walker has walked straight for a while without paying attention to the environment, and his
future movements have little dependence on his past movements.

We define the $\epsilon$-coins $(\epsilon_{i,x})_{i\in\mathbb{N}, x\in \mathbb{Z}^d}=:\epsilon$
to be iid random variables with distribution
$Q$ such that 
\[Q(\epsilon_{i,x}=1)=d\kappa \text{ and }Q(\epsilon_{i,x}=0)=1-d\kappa.\]

For fixed $\omega$, $\epsilon$, $P_{\omega,\epsilon}^x$ is the law of the Markov chain $(X_n)$ such that $X_0=x$ and that for any $e\in\mathbb{Z}^d$
such that $|e|=1$,
\[
P_{\omega,\epsilon}^x(X_{n+1}=z+e|X_n=z)
=\frac{1_{\epsilon_{n,z}=1}}{2d}+\frac{1_{\epsilon_{n,z}=0}}{1-d\kappa}[\omega(z,z+e)-\frac{\kappa}{2}].
\]
Note that the law of $X_\cdot$ under $\bar{P}_\omega^x=Q\otimes P_{\omega,\epsilon}^x$ coincides with its
law under $P_\omega^x$. Sometimes we also refer to 
$P_{\omega,\epsilon}^x (\cdot)$ as a measure on the sets of paths, without indicating the specific random path. 
Denote by
$\bar{\textbf{P}}=P\otimes Q\otimes P_{\omega,\epsilon}^o$ the law of the triple 
$(\omega, \epsilon, X_\cdot)$. 

Now we define the regeneration times in the direction $e_1$. 
Let $L$ be a fixed number which is sufficiently large.
Set $R_0=0$. Define inductively for $k\ge 0$:
\begin{align*}
&S_{k+1}=\inf\{n\ge R_k: X_{n-L}\cdot e_1>\max\{X_m\cdot e_1: m<n-L\},\\
&\qquad\qquad\qquad \epsilon_{n-i, X_{n-i}}=1, X_{n-i+1}-X_{n-i}=e_1 \text{ for all }1\le i\le L\},\\
&R_{k+1}=R\circ\theta_{S_{k+1}}+S_{k+1},
\end{align*}
where $\theta_n$ denotes the time shift of the path, i.e.,
$\theta_n X=(X_{n+i})_{i=0}^\infty$.

Let \[K=\inf\{k\ge 1: S_k<\infty,R_k=\infty\}\] and $\tau_1=\tau_1(e_1,\epsilon,X_\cdot):=S_K.$
For $k\ge 1$, the ($L$-)regeneration times are defined inductively by \[\tau_{k+1}=\tau_1\circ\theta_{\tau_k}+\tau_k .\]

By similar argument as in \cite[Lemma 2.2]{CZ1}, we can show:
\begin{lemma}\label{l7}
If $\mathbf{P}(\lim_{n\to\infty}X_n\cdot e_1/n=0)<1$, then
\begin{equation}\label{e19}
\mathbf{P}(A_{e_1}\cup A_{-e_1})=1.
\end{equation}
Moreover, on $A_{e_1}$, $\tau_i$'s are $\bar{\mathbf P}$-almost surely finite.
\end{lemma}
\pf
If $\mathbf{P}(\lim_{n\to\infty}X_n\cdot e_1/n=0)<1$, then
\[
\mathbf{P}(\varlimsup_{n\to\infty}X_n\cdot e_1/n>0)>0\quad\text{ or }\quad
\mathbf{P}(\varlimsup_{n\to\infty}X_n\cdot (-e_1)/n>0)>0.
\]
Without loss of generality, assume that 
\[\mathbf{P}(\varlimsup_{n\to\infty}X_n\cdot e_1/n>0)>0.\]
 It then follows from Theorem \ref{thm3} that $\mathbf{P}(R=\infty)>0$. 
 We want to show that 
$R_k=\infty$ for all but finitely many $k$'s.

For $k\ge 0$,
\begin{align*}
&\bar{\mathbf P}(R_{k+1}<\infty)\\
&= \bar{\mathbf P}(S_{k+1}<\infty, R\circ\theta_{S_{k+1}}<\infty)\\
&=\sum_{n,x}\bar{\mathbf P}
(S_{k+1}=n, X_n=x, R\circ\theta_n <\infty)\\
&=\sum_{n,x}E_{P\otimes Q}
\big[P_{\omega,\epsilon}(S_{k+1}=n,X_n=x)
P_{\omega,\theta^n\epsilon}^x(R<\infty)\big],
\end{align*}
where $\theta^n\epsilon$ denotes the time shift of the coins $\epsilon$, 
i.e. $(\theta^n\epsilon)_{i,x}=\epsilon_{n+i,x}$.
Note that $P_{\omega,\epsilon}(S_{k+1}=n,X_n=x)$ and 
$P_{\omega,\theta^n\epsilon}^x(R<\infty)$ are independent under the measure $Q$,
since the former is a function of $\epsilon$'s before time $n$, and the latter
involves $\epsilon$'s after time $n$.
It then follows by induction that
\begin{align*}
&\bar{\mathbf P}(R_{k+1}<\infty)\\
&=\sum_{n,x}E_P
\big[\bar{P}_\omega(S_{k+1}=n,X_n=x)
\bar{P}_\omega^x(R<\infty)\big]\\
&=\sum_{n,x}E_P\big[\bar{P}_\omega(S_{k+1}=n,X_n=x)
E_P[\bar{P}_\omega^x(R<\infty)|\omega_y: y\cdot e_1\le x\cdot e_1-L]
\big]\\
&\stackrel{(\ref{e32}), \text{Lemma }\ref{c2}}{\le} \bar{\mathbf P}(R_k<\infty)\exp{(e^{-cL})}\bar{\mathbf P}(R<\infty)\\
&\le [\exp{(e^{-cL})}\bar{\mathbf P}(R<\infty)]^{k+1},
\end{align*}
where we used in the second equality the fact that $\bar{P}_\omega(S_{k+1}=n,X_n=x)$
is $\sigma(\omega_y: y\cdot e_1\le x\cdot e_1-L)$-measurable.	
Hence, by taking $L$ sufficiently large and by the Borel-Cantelli Lemma, $\bar{\mathbf P}$-almost surely,
$R_k=\infty$ except for finitely many values of $k$.

 Let $\mathcal{O}_{e_1}$
denote the event that the signs of $X_n\cdot e_1$ change infinitely many often.
It is easily seen that (by the ellipticity of the environment)
\begin{align*}
	&\mathbf{P}(\mathcal{O}_{e_1}\cup A_{e_1}\cup A_{-e_1})=1
\shortintertext{and}
&\mathcal{O}_{e_1}\subset \{\sup_n X_n\cdot e_1=\infty\}.
\end{align*} 
However, on $\{\sup_n X_n\cdot e_1=\infty\}$, given that
$R_k$ is finite, $S_{k+1}$ is also finite.
Hence $\tau_1$ is $\bar{\mathbf P}$-almost surely finite on $\{\sup_n X_n\cdot e_1=\infty\}$,
and so are the regeneration times $\tau_2,\tau_3\ldots$. 
Therefore,
\[
\mathbf{P}(\mathcal{O}_{e_1})=\bar{\mathbf P}(\mathcal{O}_{e_1}\cap\{\tau_1<\infty\}).
\]
Since $\mathcal{O}_{e_1}\cap\{\tau_1<\infty\}=\emptyset$, we get $\mathbf{P}(\mathcal{O}_{e_1})=0$.	This gives (\ref{e19}). \qed\\

When $\mathbf{P}(R=\infty)>0$, we let 
\[
\hat{\mathbf P}(\cdot):=\bar{\mathbf P}(\cdot|R=\infty).
\]
The following proposition is a consequence of Lemma \ref{c2}.
\begin{proposition}
Assume $\mathbf{P}(R=\infty)>0$.
Let $l>r$ and
$\Lambda\subset\{x:x\cdot e_1<-r\}.$
Then for any $A\subset\Lambda\cap\{x: x\cdot e_1<-l\}$ and $k\in\mathbb{N}$,
\begin{equation}\label{prop1}
\exp(-Ce^{-\gamma l})
\le 
\dfrac{E_P\big[\bar{P}_\omega\big((X_i)_{i=0}^{\tau_k}\in\cdot, R=\infty\big)|\omega_y:y\in\Lambda\setminus A]}
{E_P\big[\bar{P}_\omega\big((X_i)_{i=0}^{\tau_k}\in\cdot, R=\infty\big)|\omega_y:y\in\Lambda]}
\le 
\exp(Ce^{-\gamma l}).
\end{equation}
Furthermore, for any $k\in\mathbb{N}$ and $n\ge 0$, $\hat{\mathbf P}$-almost surely,
\begin{equation}\label{prop2}
\exp(-e^{-cL})
\le 
\frac{\hat{\mathbf P}\big((X_{\tau_n+i}-X_{\tau_n})_{i=0}^{\tau_{n+k}-\tau_n}\in\cdot|X_{\tau_n}\big)}
{\hat{\mathbf P}\big((X_i)_{i=0}^{\tau_k}\in\cdot\big)}
\le 
\exp(e^{-cL}).
\end{equation}
\end{proposition}
\pf
First, we shall prove (\ref{prop1}).
By the definition of the regeneration times, for any finite path $x_\cdot=(x_i)_{i=0}^N, N<\infty$, there exists an event
$G_{x_\cdot}\in\sigma(\epsilon_{i,X_i},X_i: i\le N)$
 such that $G_{x_\cdot}\subset\{R>N\}$  and
\[
\{(X_i)_{i=0}^{\tau_k}=(x_i)_{i=0}^N, R=\infty\}
=G_{x_\cdot}\cap\{R\circ\theta_N=\infty\}.
\]
(For example, when $k=1$, we let 
\[
G_{x_\cdot}=\bigcup_{j=1}^\infty  
\{(X_i)_{i=0}^N=(x_i)_{i=0}^N, S_j=N, R>N\}.
\]
Then $\{(X_i)_{i=0}^{\tau_1}=(x_i)_{i=0}^N, R=\infty\}
=G_{x_\cdot}\cap\{R\circ\theta_N=\infty\}.$)

For $n\in\mathbb{N}$, we let 
\[
E_n:=G_{x_\cdot}\cap\{R\circ\theta_N\ge n\}.
\]
Note that $E_n\in\sigma(\epsilon_{i,X_i},X_i:i\le N+n)$ can be interpreted (in the sense of (\ref{2e*}))
as a set of paths with lengths $\le N+n$. Also note that $E_n\subset\{R>N+n\}$.
Then by Lemma \ref{c2} and (\ref{e32}), we have
\[
\exp(-Ce^{-\gamma l})
\le 
\dfrac{E_P\big[\bar{P}_\omega\big(E_n)|\omega_y:y\in\Lambda\setminus A]}
{E_P\big[\bar{P}_\omega\big(E_n\big)|\omega_y:y\in\Lambda]}
\le 
\exp(Ce^{-\gamma l}).
\]
(\ref{prop1}) follows by letting $n\to\infty$.

Next, we shall prove (\ref{prop2}).
Let $x\in\mathbb{Z}^d$ be any point that satisfies 
\[
\bar{\mathbf P}(X_{\tau_n}=x)>0.
\]
By the definition of the regeneration times, 
for any $m\in\mathbb{N}$,
there exists an event
$G_m^x\in\sigma\{\epsilon_{i,X_i},X_i: i\le m\}$
such that $\bar{P}_\omega(G_m^x)$ is 
$\sigma(\omega_y:y\cdot e_1\le x\cdot e_1-L)$-measurable, and
\[
\{\tau_n=m,X_m=x, R=\infty\}=G^x_m\cap\{R\circ\theta_m=\infty\}.
\]
Thus 
\begin{align}\label{2e2}
&\bar{\mathbf P}\big((X_{\tau_n+i}-X_{\tau_n})_{i=0}^{\tau_{n+k}-\tau_n}\in\cdot,X_{\tau_n}=x, R=\infty\big)\nonumber\\
&=\sum_{m}\bar{\mathbf P}\big((X_{\tau_n+i}-X_{\tau_n})_{i=0}^{\tau_{n+k}-\tau_n}\in\cdot, \tau_n=m,X_m=x,R=\infty\big)\nonumber\\
&=
\sum_{m}E_P\big[\bar{P}_\omega(G_m^x)
\bar{P}_\omega^x((X_i-x)_{i=0}^{\tau_k}\in\cdot,R=\infty)\big]\nonumber\\
&\stackrel{(\ref{prop1})}{\le}
\exp(Ce^{-\gamma L})\sum_{m}\bar{\mathbf P}(G_m^x)
\bar{\mathbf P}\big((X_i)_{i=0}^{\tau_k}\in\cdot, R=\infty\big).
\end{align}
On the other hand, 
\begin{align}\label{2e11}
\bar{\mathbf P}(X_{\tau_n}=x, R=\infty)
&=\sum_{m} E_P[\bar{P}_{\omega}(G^x_m)\bar{P}_\omega^x(R=\infty)]\nonumber\\
&\stackrel{(\ref{prop1})}{\ge}\exp(-Ce^{-\gamma L})
\sum_{m}
\bar{\mathbf P}(G_m^x)\bar{\mathbf P}(R=\infty).
\end{align}
By (\ref{2e2}) and (\ref{2e11}), we have (note that $L$ is sufficiently big)
\[
\hat{\mathbf P}\big((X_{\tau_n+i}-X_{\tau_n})_{i=0}^{\tau_{n+k}-\tau_n}\in\cdot|X_{\tau_n}=x\big)
\le 
\exp(e^{-cL})\hat{\mathbf P}\big((X_i)_{i=0}^{\tau_k}\in\cdot\big).
\]
The right side of (\ref{prop2}) is proved. The left side of (\ref{prop2}) follows likewise.
\qed

The next lemma describes the dependency of a regeneration on its remote past.
It is a version of Lemma 2.2 in \cite{CZ2}. (The denominator is omitted in the 
last equality in \cite[page 101]{CZ2}, which is corrected here, see the equality in (\ref{*8}).)

Set $\tau_0=0$.
Denote the truncated path between $\tau_{n-1}$ and $\tau_n-L$ by 
\[
P_n=(P_n^i)_{0\le i\le \tau_{n}-\tau_{n-1}-L}:=(X_{i+\tau_{n-1}}-X_{\tau_{n-1}})_{0\le i\le \tau_n-\tau_{n-1}-L}.
\]
Set
\begin{align*}
 W_n &=(\omega_{x+X_{\tau_{n-1}}})_{x\in P_n}=:\omega_{X_{\tau_{n-1}}+P_n},\\
 F_n &=X_{\tau_n}-X_{\tau_{n-1}},\\
J_n &=(P_n,W_n,F_n,\tau_n-\tau_{n-1}).
\end{align*}
For $i\ge 0$, let $h_{i+1}(\cdot|j_i,\ldots,j_1):=\hat{\mathbf{P}}(J_{i+1}\in\cdot|J_{i},\ldots,J_1)|_{J_{i}=j_i,\ldots,J_1=j_1}$ denote the transition kernel of $(J_n)$. Note that when $i=0$, $h_{i+1}(\cdot|j_i,\ldots,j_1)=h_1(\cdot|\emptyset)=\hat{\mathbf P}(J_1\in\cdot)$.
\begin{lemma}\label{l4} 
Assume $\mathbf{P}(R=\infty)>0$, $0\le k\le n$. Then $\hat{\mathbf P}$-almost surely,
\begin{equation}\label{e26}
\exp{(-e^{-c(k+1)L})}
\le 
\frac{h_{n+1}(\cdot|J_n,\ldots,J_1)}{h_{k+1}(\cdot|J_n,\ldots,J_{n-k+1})}
\le \exp{(e^{-c(k+1)L})}.
\end{equation}
\end{lemma}

\pf 
For $j_m=(p_m,w_m,f_m,t_m),m=1,\ldots n$, let
\begin{align*}
&\bar{x}_m:=f_1+\cdots+f_m,\\
&\bar{t}_m:=t_1+\cdots+t_m,\\
&B_{p_1,\ldots,p_m}:=\{R=\infty, P_i=p_i\text{ for all }i=1,\ldots,m\},\\
\text{ and }\quad
& \omega_{p_1,\ldots,p_m}:=(\omega_{\bar{x}_{i-1}+p_i})_{i=1}^m.
\end{align*}

First, we will show that for any $1\le k\le n$,
\begin{equation}\label{*8}
h_{k+1}(\cdot|j_k,\ldots,j_1)
=\frac{E_P\big[\bar{P}_\omega^{\bar{x}_k}(J_1\in \cdot,R=\infty)|\omega_{p_1,\ldots,p_k}\big]}
{E_P\big[\bar{P}_\omega^{\bar{x}_k}(R=\infty)|\omega_{p_1,\ldots,p_k}\big]}
\Big|_{\omega_{p_1,\ldots,p_k}=(w_i)_{i=1}^k}.
\end{equation}
By the definition of the regeneration times, there exists an event
\[
G_{p_1,\ldots,p_k}\in\sigma(X_{i+1},\epsilon_{i,X_i}, 0\le i\le \bar{t}_k-1)
\]
such that 
\begin{equation}\label{*7}
B_{p_1,\ldots,p_k}=G_{p_1,\ldots,p_k}\cap\{R\circ\theta_{\bar{t}_k}=\infty\}.
\end{equation}
On the one hand, for any $\sigma(J_k,\ldots, J_1)$-measurable function $g(J_k,\ldots,J_1)$,
\begin{align}\label{*15}
&E_{\bar{\mathbf P}}\big[h_{k+1}(\cdot|J_k,\ldots,J_1)g(J_k,\ldots, J_1)1_{B_{p_1,\ldots,p_k}}\big]\nonumber\\
&=E_{\bar{\mathbf P}}\big[g1_{B_{p_1,\ldots,p_k}}1_{J_{k+1}\in \cdot}\big]\nonumber\\
&=
E_P [g1_{B_{p_1,\ldots,p_k}}\bar{P}_\omega(J_{k+1}\in\cdot,B_{p_1,\ldots,p_k})]\nonumber\\
&\stackrel{(\ref{*7})}{=}
E_P \big[g1_{B_{p_1,\ldots,p_k}}\bar{P}_\omega(G_{p_1,\ldots,p_k})\bar{P}_\omega^{\bar{x}_k}(J_1\in \cdot,R=\infty)\big].
\end{align}
On the other hand, we also have
\begin{align}\label{*16}
& E_{\bar{\mathbf P}}\big[h_{k+1}(\cdot|J_k,\ldots,J_1)
g(J_k,\ldots, J_1)1_{B_{p_1,\ldots,p_k}}\big]\nonumber\\
&=
E_P \big[
h_{k+1}(\cdot|J_k,\ldots,J_1)g1_{B_{p_1,\ldots,p_k}}
\bar{P}_\omega(B_{p_1,\ldots,p_k})\big]\nonumber\\
&\stackrel{(\ref{*7})}{=}
E_P \big[
h_{k+1}(\cdot|J_k,\ldots,J_1)g1_{B_{p_1,\ldots,p_k}}
\bar{P}_\omega(G_{p_1,\ldots,p_k})\bar{P}_\omega^{\bar{x}_k}(R=\infty)\big].
\end{align}
Comparing (\ref{*15}) and (\ref{*16}) and observing that on $B_{p_1,\ldots,p_k}$, 
$\bar{P}_\omega(G_{p_1,\ldots,p_k})$ and all functions of $J_1,\ldots,J_k$
are $\sigma(\omega_y: y\in\bar{x}_{i-1}+p_i, i\le k)$-measurable
, we obtain that on $B_{p_1,\ldots,p_k}$, $P$-almost surely,
\begin{equation*}
h_{k+1}(\cdot|J_k,\ldots,J_1)
=\frac{E_P\big[\bar{P}_\omega^{\bar{x}_k}(J_1\in \cdot,R=\infty)|\omega_{\bar{x}_{i-1}+p_i},i\le k\big]}
{E_P\big[\bar{P}_\omega^{\bar{x}_k}(R=\infty)|\omega_{\bar{x}_{i-1}+p_i},i\le k\big]}.
\end{equation*}
Noting that
\[
B_{p_1,\ldots,p_k}\cap\{\omega_{p_1,\ldots,p_k}=(w_i)_{i=1}^k\}
=
\{J_i=j_i,1\le i\le k\},
\]
(\ref{*8}) is proved.

Next, we will prove the lower bound in \eqref{e26}.

When $n\ge k\ge 1$,	by formula (\ref{*8}) and (\ref{prop1}), we have
\begin{align}\label{3e2}
&h_{n+1}(\cdot| j_n,\ldots,j_1)\nonumber\\
&=
\frac{E_P[\bar{P}_\omega^{\bar{x}_n}(J_1\in\cdot,R=\infty)|\omega_{p_1,\ldots,p_n}]}
{E_P\big[\bar{P}_\omega^{\bar{x}_n}(R=\infty)|\omega_{p_1,\ldots,p_n}\big]}\bigg|_{\omega_{p_1,\ldots,p_n}=(w_i)_{i=0}^n}\nonumber\\
&\le
\frac{\exp(Ce^{-\gamma(k+1)L})
E_P[\bar{P}_\omega^{\bar{x}_n}(J_1\in\cdot,R=\infty)|\omega_{\bar{x}_{i-1}+p_i},n-k+1\le i\le n]}
{\exp(-Ce^{-\gamma(k+1)L})
E_P[\bar{P}_\omega^{\bar{x}_n}(R=\infty)|\omega_{\bar{x}_{i-1}+p_i},n-k+1\le i\le n]}\nonumber\\
&\qquad
\big|_{\omega_{p_1,\ldots,p_n}=(w_i)_{i=0}^n}\nonumber\\
&=
\exp(2Ce^{-\gamma(k+1)L})
\frac{E_P[\bar{P}_\omega^{\bar{x}_n-\bar{x}_{n-k}}(J_1\in\cdot,R=\infty)|\omega_{p_{n-k+1},\ldots, p_n}]}
{E_P[\bar{P}_\omega^{\bar{x}_n-\bar{x}_{n-k}}(R=\infty)|\omega_{p_{n-k+1},\ldots, p_n}]}\nonumber\\
&\qquad
\big|_{\omega_{p_{n-k+1},\ldots,p_n}=(w_i)_{i=n-k+1}^n}\nonumber\\
&\stackrel{(\ref{*8})}{=}
\exp(2Ce^{-\gamma(k+1)L})h_{k+1}(\cdot|j_{n},\ldots,j_{n-k+1}),
\end{align}
where we used the translation invariance of the measure $P$ in the last but one equality.

When $k=0$ and $n\ge 1$, by formula \eqref{*8} and \eqref{prop1},
\begin{align}\label{3e1}
	h_{n+1}(\cdot| j_n,\ldots,j_1)
	&\le 
	\frac{\exp(Ce^{-\gamma L})E_P[\bar{P}_\omega^{\bar{x}_n}(J_1\in\cdot,R=\infty)]}
	{\exp(-Ce^{-\gamma L})E_P[\bar{P}_\omega^{\bar{x}_n}(R=\infty)]}\nonumber\\
	&=\exp(2Ce^{-\gamma L})\hat{\mathbf P}(J_1\in\cdot)\nonumber\\
	&=\exp(2Ce^{-\gamma L})h_1(\cdot|\emptyset).
\end{align}

When $k=n=0$, \eqref{e26} is trivial.
Hence combining \eqref{3e2} and \eqref{3e1}, the lower bound in (\ref{e26}) follows as we take $L$ sufficiently big. The upper bound follows likewise.\qed
 
\begin{lemma}\label{l6}
Suppose that a sequence of non-negative random variables $(X_n)$ satisfies
\[
a\le \frac{\ud P(X_{n+1}\in\cdot|X_1,\ldots, X_n)}{\ud \mu}\le b
\]
for all $n\ge 1$, where $a\le 1\le b$ are constants and $\mu$ is a probability measure. Let $m_\mu\le \infty$ be
the mean of $\mu$. Then almost surely,
\begin{equation}\label{e21}
a m_\mu\le\varliminf_{n\to\infty}
\frac{1}{n}\sum_{i=1}^n X_i
\le
\varlimsup_{n\to\infty}
\frac{1}{n}\sum_{i=1}^n X_i
\le b m_\mu.
\end{equation}
\end{lemma}

Before giving the proof, let us recall the ``splitting representation" of random variables:
\begin{proposition}\cite[Page 94]{Tho}\label{prop}
Let $\nu$ and $\mu$ be probability measures. Let $X$ be a random variable with law $\nu$. If 
for some $a\in(0,1)$,
\[
\frac{\ud\nu}{\ud\mu}\ge a,
\]
then, enlarging the probability space if necessary, we can find independent
random variables $\Delta, \pi, Z$ such that
\begin{itemize}
\item[i)] $\Delta$ is Bernoulli with parameter $1-a$, i.e., $P(\Delta=1)=1-a$,
$P(\Delta=0)=a$;
\item[ii)] $\pi$ is of law $\mu$, and $Z$ is of law $(\nu-a\mu)/(1-a)$;
\item[iii)] $X=(1-\Delta)\pi+\Delta Z$.
\end{itemize}
\end{proposition}

\noindent{\it Proof of Lemma \ref{l6}}:\\
By Proposition \ref{prop}, enlarging the probability space if necessary, there are
random variables $\Delta_i,\pi_i,Z_i,i\ge 1$,
such that for any $i\in\mathbb{N}$,
\begin{itemize}
\item 
$\Delta_i$ is Bernoulli with parameter $(1-a)$, and $\pi_i$ is of law $\mu$;
\item 
$\Delta_i, \pi_i$ and $Z_i$ are mutually independent;
\item 
$(\Delta_i,\pi_i)$
is independent of
$\sigma(\Delta_k,\pi_k,Z_k: k<i)$;
\item
$X_i=(1-\Delta_i)\pi_i+\Delta_i Z_i$.
\end{itemize}
Note that since $X_i$'s are supported on $[0,\infty)$, $\pi_i\ge 0$ and $Z_i\ge 0$ for all $i\in\mathbb{N}$.
Thus by the law of large numbers, almost surely,
\[
\varliminf_{n\to\infty}\frac{1}{n}\sum_{i=1}^n X_i\ge
\lim_{n\to\infty}\frac{1}{n}\sum_{i=1}^n (1-\Delta_i)\pi_i=a m_\mu.
\]
This proves the first inequality of (\ref{e21}). 

If $m_\mu=\infty$, the last inequality
of (\ref{e21}) is trivial. Assume that $m_\mu<\infty$.
Let $(\tilde{\Delta}_i)_{i\ge 1}$ be an iid Bernoulli sequence with
parameter $1-b^{-1}$ such that every $\tilde{\Delta}_i$ is independent of
all the $X_n$'s.
By a similar splitting procedure, we can construct non-negative
random variables $\tilde{\pi}_i,\tilde{Z}_i, i\ge 1$,
such that $(\tilde{\pi}_i)_{i\ge 1}$ are iid with law $\mu$, and
\[
\tilde{\pi}_i=(1-\tilde{\Delta}_i)X_i+\tilde{\Delta}_i\tilde{Z}_i.
\]

Let $Y_i=(1-b^{-1}-\tilde{\Delta}_i)X_i 1_{X_i\le i}$, we will first show that
\begin{equation}\label{e22}
\lim_{n\to\infty}\frac{1}{n}\sum_{i=1}^n Y_i=0.
\end{equation}
By Kronecker's Lemma, it suffices to show that
\[
\sum_{i=1}^\infty \frac{Y_i}{i} \text{ converges.}
\]
Observe that $(\sum_{i=1}^n Y_i/i)_{n\in\mathbb{N}}$ is a martingale sequence.
Moreover, for all $n\in\mathbb{N}$,
\begin{align*}
E\big(\sum_{i=1}^n \frac{Y_i}{i}\big)^2
= \sum_{i=1}^n EY_i^2/i^2
&\le \sum_{i=1}^\infty EX_i^2 1_{X_i\le i}/i^2\\
&\le b \sum_{i=1}^\infty E\tilde{\pi}_i^2 1_{\tilde{\pi}_i\le i}/i^2\\
&= b\int_0^\infty x^2 (\sum_{i\ge x}\frac{1}{i^2})\ud\mu\\
&\le C\int_0^\infty x\ud\mu=Cm_\mu<\infty.
\end{align*}
By the $L^2$-martingale convergence theorem, $\sum Y_i/i$
converges a.s. and in $L^2$. This proves (\ref{e22}). 

Since
\[
\sum_i P(Y_i\neq (1-b^{-1}-\tilde{\Delta}_i)X_i)
\le 
\sum_i P(X_i>i)
\le 
b \sum_i P(\pi_1>i)\le b m_\mu<\infty,
\]
by the Borel-Cantelli lemma, it follows from (\ref{e22}) that
\[
\lim_{n\to\infty}\frac{1}{n}\sum_{i=1}^n (1-b^{-1}-\tilde{\Delta}_i)X_i=0, \text{a.s.}.
\]
Hence almost surely,
\begin{equation*}
m_\mu=\lim_{n\to\infty}\frac{1}{n}\sum_{i=1}^n \tilde{\pi}_i
\ge 
\varlimsup_{n\to\infty}\frac{1}{n}\sum_{i=1}^n (1-\tilde{\Delta}_i)X_i
=\varlimsup_{n\to\infty}\frac{1}{n}\sum_{i=1}^n b^{-1} X_i.
\end{equation*}
The last inequality of (\ref{e21}) is proved.\qed 

\begin{theorem}\label{lln}
There exist two deterministic numbers $v_{e_1},v_{-e_1}\ge 0$ such that $\mathbf{P}$-almost surely,
\begin{equation}\label{e25}
\lim_{n\to\infty}\frac{X_n\cdot e_1}{n}=v_{e_1} 1_{A_{e_1}}-v_{-e_1}1_{A_{-e_1}}.
\end{equation}
Moreover, if $v_{e_1}>0$, then $E_{\hat{\mathbf P}}\tau_1<\infty$ and 
$\mathbf{P}(A_{e_1}\cup A_{-e_1})=1$.
\end{theorem}
\pf
We only consider the nontrivial case that $\mathbf{P}(\lim X_n\cdot e_1/n=0)<1$,
which by Lemma \ref{l7} implies
$\mathbf{P}(A_{e_1}\cup A_{-e_1})=1$. 
Without loss of generality, assume $\mathbf{P}(\varlimsup_{n\to\infty}X_n\cdot e_1/n>0)>0$.
We will show that on $A_{e_1}$,
\[
\lim_{n\to\infty}X_n\cdot e_1/n=v_{e_1}>0, \text{ $\mathbf{P}$-a.s..}
\]

By (\ref{prop2}) and Lemma \ref{l6}, we obtain that
$\mathbf{P}(\cdot|A_{e_1})$-almost surely,
\begin{align}
\exp{(-e^{-cL})}E_{\hat{\mathbf P}}X_{\tau_1}\cdot e_1
&\le \varliminf_{n\to\infty}\frac{X_{\tau_n}\cdot e_1}{n}\nonumber\\
&\le \varlimsup_{n\to\infty}\frac{X_{\tau_n}\cdot e_1}{n}
\le \exp{(e^{-cL})}E_{\hat{\mathbf P}}X_{\tau_1}\cdot e_1,\label{e23}\\
\exp{(-e^{-cL})}E_{\hat{\mathbf P}}\tau_1
&\le \varliminf_{n\to\infty}\frac{\tau_n}{n}
\le \varlimsup_{n\to\infty}\frac{\tau_n}{n}
\le \exp{(e^{-cL})}E_{\hat{\mathbf P}}\tau_1. \label{e24}
\end{align}
Note that (\ref{e23}), (\ref{e24}) hold even if 
$E_{\hat{\mathbf P}}X_{\tau_1}\cdot e_1=\infty$ or $E_{\hat{\mathbf P}}\tau_1=\infty$.
But it will be shown later that under our assumption, both of them are finite.

We claim that 
\begin{equation}\label{e5}
E_{\hat{\mathbf P}}X_{\tau_1}\cdot e_1<\infty.
\end{equation}
To see this, let $\Theta:=\{i: X_{\tau_k}\cdot e_1=i \text{ for some }k\in\mathbb{N}\}$.
Since $\tau_i$'s are finite on $A_{e_1}$,
there exist (recall that $\tau_0=0$) a sequence $(k_n)_{n\in\mathbb{N}}$ such that
$X_{\tau_{k_n}}\cdot e_1\le n<X_{\tau_{k_n+1}}\cdot e_1$ for all $n\in\mathbb{N}$ and 
$\lim_{n\to\infty}k_n=\infty$.
Hence for $n\ge 1$,
\[
\frac{\sum_{i=1}^n 1_{i\in \Theta}}{n}\le 
\frac{k_n+1}{X_{\tau_{k_n}}\cdot e_1}, \quad\text{ $\hat{\mathbf P}$-a.s..}
\]
Then, $\hat{\mathbf P}$-a.s.,
\[
\varlimsup_{n\to\infty}
\frac{\sum_{i=1}^n 1_{i\in \Theta}}{n}
\le
\varlimsup_{n\to\infty}\frac{n}{X_{\tau_n}\cdot e_1}.
\]

Let $B_k=\{\epsilon_{k,X_k}=0, X_{k+1}-X_k=e_1, \epsilon_{k+i,X_{k+i}}=1,X_{k+i+1}-X_{k+i}=e_1
 \text{ for all }1\le i\le L\}$. Then
\[
 \bar{P}_\omega(B_k)
 \ge 
 (d\kappa)^L(1-d\kappa)(\frac{\kappa}{2})(\frac{1}{2d})^L
 \stackrel{1\ge 2d\kappa}{>}(\frac{\kappa}{2})^{L+2}.
\]
Observe that by the definition of the regeneration times, for $n> L+1$,
\begin{align*}
&\{T_{n-L-1}=k,X_k= x-(L+1)e_1, R>k\}\cap B_k\cap\{R\circ\theta_{k+L+1}=\infty\}\\
&\subset\{R=\infty, n\in\Theta, T_n=k+L+1,X_{T_n}=x\}.
\end{align*}
Hence for $n> L+1$,
\begin{align*}
& \hat{\mathbf P}(n\in \Theta)\\
&\ge  
\sum_{k\in\mathbb{N},x\in\mathcal{H}_n}
\hat{\mathbf P}(B_k\cap\{T_{n-L-1}=k,X_k= x-(L+1)e_1,R\circ\theta_{k+L+1}=\infty\})\\
&\ge  \sum_{k\in\mathbb{N},x\in\mathcal{H}_n}
E_P \big[P_\omega\big(T_{n-L-1}=k,X_k= x-(L+1)e_1,R>k\big)(\frac{\kappa}{2})^{L+2}\\
&\qquad\qquad\qquad\qquad\qquad\qquad\qquad\quad\times 
P_\omega^x(R=\infty)\big]/\mathbf{P}(R=\infty).
\end{align*}
Since by (\ref{prop1}) and the translation invariance of $P$,
\[
E_P\big[P_\omega^x(R=\infty)|\omega_y:y\cdot e_1\le x\cdot e_1-L-1\big]
\ge 
\exp(-e^{-cL})\mathbf{P}(R=\infty),
\]
we have for $n>L+1$,
\begin{align}\label{2e10}
& \hat{\mathbf P}(n\in \Theta)\nonumber\\
&\ge
(\frac{\kappa}{2})^{L+2}\exp(-e^{-cL})
\sum_{k\in\mathbb{N},x\in\mathcal{H}_n}\mathbf{P}(T_{n-L-1}=k,X_k= x-(L+1)e_1,R>k)\nonumber\\
&\ge (\frac{\kappa}{2})^{L+2} e^{-1}\mathbf{P}(R=\infty).
\end{align}
Hence
\begin{align*}
\frac{C}{E_{\hat{\mathbf P}}X_{\tau_1}\cdot e_1}
\stackrel{(\ref{e23})}{\ge} E_{\hat{\mathbf P}}\varlimsup_{n\to\infty}\frac{n}{X_{\tau_n}\cdot e_1}
&\ge  E_{\hat{\mathbf P}}\varlimsup_{n\to\infty}
\frac{\sum_{i=1}^n 1_{i\in \Theta}}{n}\\
&\ge \varlimsup_{n\to\infty}E_{\hat{\mathbf P}}\frac{\sum_{i=1}^n 1_{i\in \Theta}}{n}\\
&\stackrel{(\ref{2e10})}{\ge} (\frac{\kappa}{2})^{L+2} e^{-1}\mathbf{P}(R=\infty)>0.
\end{align*}
This gives (\ref{e5}).

Now we can prove the theorem.
By (\ref{e23}) and (\ref{e24}), 
\begin{align}\label{*9}
\exp{(-2e^{-cL})}\frac{E_{\hat{\mathbf P}}X_{\tau_1}\cdot e_1}{E_{\hat{\mathbf P}}\tau_1}
&\le \varliminf_{n\to\infty}\frac{X_{\tau_n}\cdot e_1}{\tau_{n+1}}\nonumber\\
&\le \varlimsup_{n\to\infty}\frac{X_{\tau_{n+1}}\cdot e_1}{\tau_n}
\le \exp{(2e^{-cL})}\frac{E_{\hat{\mathbf P}}X_{\tau_1}\cdot e_1}{E_{\hat{\mathbf P}}\tau_1},
\end{align}
$\mathbf{P}(\cdot|A_{e_1})$-almost surely.
Further, by the fact that $|X_i|\le i$ and the obvious inequalities
\begin{equation*}
\varliminf_{n\to\infty}\frac{X_{\tau_n}\cdot e_1}{\tau_{n+1}}
\le 
\varliminf_{n\to\infty}\frac{X_n\cdot e_1}{n}
\le 
\varlimsup_{n\to\infty}\frac{X_n\cdot e_1}{n}
\le 
\varlimsup_{n\to\infty}\frac{X_{\tau_{n+1}}\cdot e_1}{\tau_n},
\end{equation*}
we have that
\begin{equation*}
\varlimsup_{n\to\infty}
\Bigl\lvert \frac{X_n\cdot e_1}{n}-
\frac{E_{\hat{\mathbf P}}X_{\tau_1}\cdot e_1}{E_{\hat{\mathbf P}}\tau_1}\Bigr\rvert
\le \exp{(2e^{-cL})}-1, \text{ $\mathbf{P}(\cdot|A_{e_1})$-a.s.}
\end{equation*}
Therefore, $\mathbf{P}(\cdot|A_{e_1})$-almost surely,
\[
\lim_{n\to\infty} \frac{X_n\cdot e_1}{n}
=
\lim_{L\to\infty}
\frac{E_{\hat{\mathbf P}}X_{\tau_1^{(L)}}\cdot e_1}{E_{\hat{\mathbf P}}\tau_1^{(L)}}:=v_{e_1},
\]
where $\tau_1$ is written as $\tau_1^{(L)}$ to indicate that it is an $L$-regeneration time.
Moreover, our assumption $\mathbf{P}(\varlimsup_{n\to\infty}X_n\cdot e_1/n>0)>0$ implies that
$v_{e_1}>0$ and (by (\ref{*9}))
\[
E_{\hat{\mathbf P}}\tau_1<\infty.
\]
Our proof is complete.\qed\\

If $v_{e_1}>0$, then it follows by (\ref{e24}) that
\begin{equation}\label{etau}
E_{\hat {\mathbf P}}\tau_n\le CnE_{\hat{\mathbf P}}\tau_1<\infty.
\end{equation}

Observe that although Theorem \ref{lln} is stated for $e_1$, the previous arguments, if properly modified, still work if one replaces
$e_1$ with any $z\in\mathbb{R}^d\setminus\{o\}$. So Theorem \ref{lln} is true for the general case. That is, for any $z\neq o$, there exist two deterministic constants $v_z, v_{-z}\ge 0$ such that
\[
\lim_{n\to\infty}\frac{X_n\cdot z}{n}=v_z1_{A_z}-v_{-z}1_{A_{-z}} 
\]
and that $\mathbf{P}(A_z\cup A_{-z})=1$ if $v_z>0$.
Then, by the same argument as in \cite[page 1112]{Go}, one concludes that the
limiting velocity $\lim_{n\to\infty}X_n/n$ can take at most two antipodal values.
This proves display \eqref{LVclln} of Theorem \ref{thm2}.

\section{Heat kernel estimate}\label{sechke}
The following heat kernel estimates are crucial for the proof of the uniqueness of the
non-zero velocity in the next section. Although in the mixing case we don't have iid
regeneration slabs, we know that (by Lemma \ref{l4}) a regeneration slab has little dependence on its remote
past. This allows us to use coupling techniques to get the same heat
kernel estimates as in \cite{Be}:
\begin{theorem}[Heat kernel estimate]\label{hke}
Assume $v_{e_1}>0$. For $x\in\mathbb{Z}^d$ and $n\in \mathbb{N}$,
we let 
\[Q(n,x):=\hat{\mathbf P}(x \text{ is visited in }[\tau_{n-1},\tau_n)).\]
Then, for any $x\in\mathbb{Z}^d$ and $n\in \mathbb{N}$,
\begin{align}
&\hat{\mathbf P}(X_{\tau_n}=x)\le Cn^{-d/2},\label{ehke}\\
&
\sum_{x\in\mathbb{Z}^d}Q(n,x)^2\le C(E_{\hat{\mathbf P}}\tau_1)^2 n^{-d/2}.\label{ehke2}
\end{align}
\end{theorem}

By Lemma \ref{l4}, we have for $n\ge 2$ and $1\le k\le n-1$, $\hat{\mathbf P}$-almost surely,
\begin{align}\label{new1}
\frac{h_{k+1}(\cdot|J_{n-1},\ldots,J_{n-k})}{h_k(\cdot|J_{n-1},\ldots,J_{n-k+1})}
&=
\frac{h_{k+1}(\cdot|J_{n-1},\ldots,J_{n-k})}{h_n(\cdot|J_{n-1},\ldots,J_1)}
\frac{h_n(\cdot|J_{n-1},\ldots,J_1)}{h_k(\cdot|J_{n-1},\ldots,J_{n-k+1})}\nonumber\\
&\ge\exp(-e^{-c(k+1)L}-e^{-ckL})\nonumber\\
&\ge 1-e^{-ckL}
\end{align}
for large $L$. Hence for $n\ge 2$ and $1\le k\le n-1$, we can define a (random) probability measure $\zeta_{n,k}^{J_{n-1},\ldots,J_{n-k}}$ that satisfies
\begin{align}\label{new2}
\MoveEqLeft
h_{k+1}(\cdot|J_{n-1},\ldots,J_{n-k})\\
&=e^{-ckL}\zeta_{n,k}^{J_{n-1},\ldots,J_{n-k}}(\cdot)+(1-e^{-ckL})h_k(\cdot|J_{n-1},\ldots, J_{n-k+1}).\nonumber
\end{align}

To prove Theorem \ref{hke}, we will construct in Section~\ref{J} a sequence of
random variables $(\tilde J_i, i\in\mathbb{N})$ such that for any
$n\in\mathbb{N}$, 
\begin{equation}\label{new0}
(\tilde J_1,\ldots,\tilde J_n)\sim \hat{\mathbf P}(J_1\in\cdot,\ldots, J_n\in\cdot),
\end{equation}
where ``$X\sim\mu$" means ``$X$ is of law $\mu$". 
\subsection{Construction of the $\tilde J_i$'s}\label{J}
Our construction consists of three steps:

\noindent{\it Step 1.}
We let $\tilde J_1, \tilde J_{2,1}, \tilde \Delta_{2,1}$ be independent random variables  such that 
\[
\tilde J_1\sim h_1(\cdot|\emptyset),\quad \tilde J_{2,1}\sim h_1(\cdot|\emptyset)
\] 
and $\tilde \Delta_{2,1}$ is Bernoulli with parameter $e^{-cL}$. Let $\tilde Z_{2,1}$ be independent of $\sigma(\tilde J_{2,1}, \tilde \Delta_{2,1})$ such that
\[
P(\tilde Z_{2,1}\in\cdot|\tilde J_1)=
\zeta_{2,1}^{\tilde J_1}(\cdot).
\]
Setting $\tilde J_{2}:=(1-\tilde \Delta_{2,1})\tilde J_{2,1}+\tilde\Delta_{2,1}\tilde Z_{2,1}$, by \eqref{new2} we have
\[
(\tilde J_1, \tilde J_2)\sim
\hat{\mathbf P}(J_1\in\cdot,J_2\in\cdot).
\]
\noindent{\it Step 2.}
For $n\ge 3$, assume we have constructed $\tilde J_1$ and $(\tilde J_{i,1},
\tilde\Delta_{i,j}, \tilde Z_{i,j}, 1\le j<i\le n-1)$ such that
\[
(\tilde J_1,\ldots, \tilde J_{n-1})
\sim
\hat{\mathbf P}(J_1\in\cdot,\ldots,J_{n-1}\in\cdot),
\]
where for $2\le j\le i\le n-1$,
\[
\tilde J_{i,j}:=(1-\tilde\Delta_{i,j-1})\tilde J_{i,j-1}+\tilde\Delta_{i,j-1}\tilde Z_{i,j-1}
\]
and
\[
\tilde J_i:=
\tilde J_{i,i}.
\]
Then, we define $\tilde J_{n,1}$ and $(\tilde\Delta_{n,k}, \tilde Z_{n,k}, 1\le
k<n)$ to be random variables such that, \textit{conditioning} on the values
of $\tilde J_1$ and $(\tilde J_{i,1}, \tilde\Delta_{i,j}, \tilde Z_{i,j}, 1\le
j<i<n)$,
\begin{itemize}
\item $(\tilde J_{n,1}, \tilde\Delta_{n,k}, \tilde Z_{n,k}, 1\le k\le n-1)$ are conditionally independent;
\item The conditional distribution of $\tilde J_{n,1}$ is $h_1(\cdot|\emptyset)$;
\item For $1\le k\le n-1$, the conditional distribution of $\tilde Z_{n,k}$ is $\zeta_{n,k}^{\tilde J_{n-1},\ldots, \tilde J_{n-k}}(\cdot)$, and $\tilde\Delta_{n,k}$ is Bernoulli with parameter $e^{-ckL}$.
\end{itemize}
\noindent{\it Step 3.} For $2\le k\le n$, set 
\begin{align*}
&\tilde J_{n,k}:=(1-\tilde\Delta_{n,k-1})\tilde
J_{n,k-1}+\tilde\Delta_{n,k-1}\tilde Z_{n,k-1},\\
&\tilde J_n:=\tilde J_{n,n}.
\end{align*}
Then (by \eqref{new2}) almost surely,
\begin{equation}\label{3e3}
P(\tilde J_{n,k}\in\cdot|\tilde J_{n-1},\ldots,\tilde J_1)=h_k(\cdot|\tilde J_{n-1},\ldots,\tilde J_{n-k+1}).
\end{equation}
It follows immediately that
\begin{equation*}
(\tilde J_1,\ldots,\tilde J_n)\sim \hat{\mathbf P}(J_1\in\cdot,\ldots, J_n\in\cdot).
\end{equation*}
Therefore, by induction, we have constructed $(\tilde J_i, i\in\mathbb{N})$ such that \eqref{new0} holds for all $n\in\mathbb{N}$. 

In what follows,  with abuse of notation, we will identify $\tilde J_i$ with
$J_i$ and simply write $\tilde J_{i,j}, \tilde \Delta_{i,j}, \tilde Z_{i,j}$ as
$J_{i,j}, \Delta_{i,j}$ and $Z_{i,j}$, $1\le j<i$. We still use $\hat{\mathbf
P}$ to denote the law of the random variables in the enlarged probability space.
\begin{remark}
To summarize, we have introduced random variables $J_{i,j}, \Delta_{i,j},
Z_{i,j}$, $1\le j< i$ such that for any $n\ge 2$,
\begin{align*}
&J_{n,2}=(1-\Delta_{n,1})J_{n,1}+\Delta_{n,1}Z_{n,1},\\
&\ldots,\\
&J_{n,n-1}=(1-\Delta_{n,n-2})J_{n,n-2}+\Delta_{n,n-2}Z_{n,n-2},\\
&J_n=(1-\Delta_{n,n-1})J_{n,n-1}+\Delta_{n,n-1}Z_{n,n-1}.
\end{align*}
Intuitively, we flip a sequence of ``coins" $\Delta_{n,n-1},\ldots,\Delta_{n,1}$ to determine whether $J_1,\ldots,J_{n-1}$ are in the ``memory" of $J_n$. For instance, if 
\[
\Delta_{n,n-1}=\cdots=\Delta_{n,n-i}=0,
\]
 then $J_n=J_{n,n-i}$ doesn't ``remember" $J_1,\ldots, J_i$ (in the sense that
\[
\hat{\mathbf P}(J_{n,n-i}\in\cdot|J_{n-1},\ldots, J_1)
=
h_{n-i}(\cdot|J_{n-1},\ldots, J_{i+1}).
\]
See \eqref{3e3}.).
\end{remark}
\subsection {Proof of Theorem \ref{hke}}
For $1<i\le n$, let $I_n(i)$ be the event that 
$\Delta_{i,i-1}=\ldots=\Delta_{i,1}=0$ and
$\Delta_{m,m-1}=\ldots=\Delta_{m,m-i}=0$ for all
$i<m\le n$.
Note that on $I_n(i)$, 
\begin{equation}\label{new3}
J_i=J_{i,1}\text{ and }J_m=J_{m,m-i}
\text{ for all }i<m\le n.
\end{equation}
\begin{lemma}\label{liid}
For $n\ge 2$, let $H$ be a nonempty subset of $\{2,\ldots, n\}$, and set
\[M_n:=\{1< i< n: I_n(i)\neq\emptyset\}.\]
 Conditioning on the
event $\{M_n=H\}$, the sequence $(J_i)_{i\in H}$ is iid and independent of
$(J_i)_{i\in \{1,\ldots,n\}\setminus H}$.
\end{lemma}
\noindent{\it Proof of Lemma \ref{liid}:}
From our construction, it follows that for any $i>1$, $J_{i,1}$ is
independent of 
\[
\sigma(\Delta_{k,j}, 1\le j<k)
\vee
\sigma(J_l, 1\le l<i)
\vee
\sigma(J_{m,m-i},m>i).
\]
Hence, by \eqref{new3}, for any $i\in H$ and any appropriate measurable sets
$(V_j)_{1\le j\le n}$,
\begin{align*}
&\hat{\mathbf P}(J_j\in V_j, 1\le j\le n|M_n=H)\\
&=\hat{\mathbf P}(J_{i,1}\in V_i)
\hat{\mathbf P}(J_j\in V_j, 1\le j\le n, j\neq i|M_n=H).
\end{align*}
By induction, we get
\begin{align*}
&\hat{\mathbf P}(J_j\in V_j, 1\le j\le n|M_n=H)\\
&=\prod_{i\in H}\hat{\mathbf P}(J_{i,1}\in V_i)
\hat{\mathbf P}(J_j\in V_j, 1\le j\le n, j\notin H|M_n=H).
\end{align*}
The lemma is proved.\qed\\

\noindent{\it Proof of Theorem \ref{hke}:}
By Lemma \ref{liid}, for $i\in H\subset\{2,\ldots, n\}$,
\[
\hat{\mathbf P}\big(X_{\tau_i}-X_{\tau_{i-1}}=(L+1)e_1\pm e_j|M_n=H\big)
=
\hat{\mathbf P}(X_{\tau_1}=(L+1)e_1\pm e_j)>0
\]
for all $j\in\{1,\ldots, d\}$, where the last inequality  is due to ellipticity.
Hence arguing as in \cite[pages 736, 737]{Be}, using Lemma~\ref{liid} and the
heat kernel estimate for bounded iid random walks in $\mathbb{Z}^d$, we get that
for any $x\in\mathbb{Z}^d$,
\[
\hat{\mathbf P}(\sum_{i\in H}X_{\tau_i}-X_{\tau_{i-1}}=x|M_n=H)
\le C|H|^{-d/2},
\]
where $|H|$ is the cardinality of $H$.
Hence, for any subset $H\subset\{2,\ldots,n\}$ such that $|H|\ge n/2$,
\begin{align}\label{e13}
& \hat{\mathbf P}(X_{\tau_n}=x|M_n=H)\nonumber\\
&=\sum_y \hat{\mathbf P}\big(\sum_{i\in H}X_{\tau_i}-X_{\tau_{i-1}}=x-y, 
\sum_{i\in \{1,\ldots,n\}\setminus
H}X_{\tau_i}-X_{\tau_{i-1}}=y|M_n=H\big)\nonumber\\
&=\sum_y 
\bigg[\hat{\mathbf P}\big(\sum_{i\in H}X_{\tau_i}-X_{\tau_{i-1}}=x-y|M_n=H\big)\nonumber\\
&\qquad\qquad\qquad\qquad\qquad\times
\hat{\mathbf P}\big(\sum_{i\in \{1,\ldots,n\}\setminus
H}X_{\tau_i}-X_{\tau_{i-1}}=y|M_n=H\big)\bigg]\nonumber\\
&\le  C n^{-d/2},
\end{align}
where we used Lemma~\ref{liid} in the second equality.

On the other hand, 
\begin{align*}
|M_n|
&\ge n-\sum_{i=2}^n \bigg(
1_{\Delta_{i,i-1}+\cdots+\Delta_{i,1}>0}+\sum_{m=i+1}^n
1_{\Delta_{m,m-1}+\cdots+\Delta_{m,m-i}>0}
\bigg)\\
&=n-\sum_ {i=2} ^n
1_{\Delta_{i,i-1}+\cdots+\Delta_{i,1}>0}-\sum_{m=2}^n\sum_{i=2}^{m-1}1_{\Delta_{
m,m-1}+\cdots+\Delta_{m,m-i}>0}\\
&\ge n-2\sum_{m=2}^n K_m,
\end{align*}
where $K_m:=\sup\{1\le j<m:\Delta_{m,j}=1\}$. Here we follow the convention
that $\sup\emptyset=0$.
Since $K_m$'s are independent, and for $m\ge 2$,
\begin{align*}
E e^{K_m} &= \sum_{j=0}^{m-1} e^j \hat{\mathbf P}(K_m=j)\\
&\le  \sum_{j=1}^{m-1} e^j \hat{\mathbf P}(\Delta_{m,j}=1)+1\\
&\le \sum_{j=1}^\infty e^j e^{-cjL}+1\to 1 \text{ as $L\to\infty$},
\end{align*}
we take $L$ to be large enough such that $E e^{K_m}\le e^{1/8}$ for all
$m\ge 2$ and so
\begin{align}\label{e14}
\hat{\mathbf P}(|M_n|<n/2)
&\le \hat{\mathbf P}(K_2+\cdots+K_n>n/4)\nonumber\\
&\le e^{-n/4}E e^{K_2+\cdots +K_n}
\le e^{-n/8}.
\end{align}
By (\ref{e13}) and (\ref{e14}), inequality (\ref{ehke}) follows immediately. 

Furthermore, since
\begin{align*}
& Q(n,x)\\
&= \sum_y \hat{\mathbf P}(X_{\tau_{n-1}}=y)
\hat{\mathbf P}(x \text{ is visited in }[\tau_{n-1},\tau_n)|X_{\tau_{n-1}}=y)\\
&\stackrel{\text{Lemma }\ref{l4}}{\le} 
C\sum_y \hat{\mathbf P}(X_{\tau_{n-1}}=y)\hat{\mathbf P}((x-y)\text{ is visited during }[0,\tau_1)),
\end{align*}
by H\"{o}lder's inequality we have
\begin{align*}
& Q(n,x)^2\\
&\le  C
\big[\sum_y\hat{\mathbf P}\big((x-y)\text{ is visited during }[0,\tau_1)\big)\big]\\
&\qquad\qquad\times\big[\sum_y\hat{\mathbf P}(X_{\tau_{n-1}}=y)^2
\hat{\mathbf P}\big((x-y)\text{ is visited during }[0,\tau_1)\big)\big]\\
&\le CE_{\hat{\mathbf P}}\tau_1 \sum_y \hat{\mathbf P}(X_{\tau_{n-1}}=y)^2
\hat{\mathbf P}\big((x-y)\text{ is visited during }[0,\tau_1)\big).
\end{align*}
Hence
\begin{align*}
& \sum_x Q(n,x)^2\\
&\le  CE_{\hat{\mathbf P}}\tau_1 
\sum_y \big[\hat{\mathbf P}(X_{\tau_{n-1}}=y)^2
\sum_x\hat{\mathbf P}\big((x-y)\text{ is visited during }[0,\tau_1)\big)\big]\\
&\le C(E_{\hat{\mathbf P}}\tau_1)^2 \sum_y \hat{\mathbf P}(X_{\tau_{n-1}}=y)^2\\
&\stackrel{(\ref{ehke})}{\le}  C(E_{\hat{\mathbf P}}\tau_1)^2 n^{-d/2}\sum_y \hat{\mathbf P}(X_{\tau_{n-1}}=y)
=C(E_{\hat{\mathbf P}}\tau_1)^2 n^{-d/2}.
\end{align*}
Theorem \ref{hke} is proved.\qed

\section{The uniqueness of the non-zero velocity}\label{secunique}
In this section we will show that in high dimension ($d\ge 5$), there
exists at most one non-zero velocity. The idea is the following.
Consider two random walk paths: one starts at the origin, the other starts near 
the $n$-th regeneration position of the first path. By Levy's martingale convergence
theorem, the second path is ``more and more transient" as $n$ grows (Lemma \ref{l1}).
On the other hand, by heat kernel estimates, when $d\ge 5$, two ballistic walks in opposite directions 
will grow further and further apart from each other (see Lemma \ref{l3}), thus they are almost independent.
This contradicts the previous fact that starting at the $n$-th regeneration point of the first path will prevent the second path from being transient in the opposite direction.\\

Set
$\delta=\delta(d):=\frac{d-4}{8(d-1)}$ (the reason of choosing this notation
will become clear in (\ref{2e6}).).
For any finite path $y_\cdot=(y_i)_{i=0}^M, M<\infty$, define $A(y_\cdot, z)$ to be the
set of paths $(x_i)_{i=0}^N, N\le\infty$ that satisfy
\begin{itemize}
\item[1)] $x_0=y_0+z$;
\item[2)] $d(x_i,y_j)>(i\vee j)^\delta$ if $i\vee j>|z|/3$.
\end{itemize}

The motivation for the definition of $A(y_\cdot, z)$ is as follows.
Note that for two paths $x_\cdot=(x_i)_{i=0}^N$ and $y_\cdot=(y_i)_{i=0}^M$ with $x_0=y_0+z$, if $i\vee j\le |z|/3$, then
\[
d(x_i,y_j)
\ge d(x_0,y_0)-d(x_0,x_i)-d(y_0,y_j)
\ge |z|-i-j
\ge |z|/3.
\]
Hence, for $(x_i)_{i=0}^N\in A(y_\cdot, z)$,
\begin{align}\label{e10}
\sum_{i\le N,j\le M}e^{-\gamma d(x_i,y_j)}
&\le 
\sum_{0\le i,j\le |z|/3}e^{-\gamma |z|/3}+\sum_{i\vee j>|z|/3}e^{-\gamma(i^\delta+j^\delta)/2}
\nonumber\\
&\le 
(\frac{|z|}{3})^2 e^{-\gamma |z|/3}+(\sum_{i=0}^\infty e^{-\gamma i^\delta/2})^2<C.
\end{align}
This gives us (by ($*$)) an estimate of the interdependence between
$\sigma\big(\omega_x: x\in (x_i)_{i=0}^N\big)$ and $\sigma\big(\omega_x: x\in (y_i)_{i=0}^M\big)$.

In what follows, we use 
\[
\tau'_\cdot=\tau_\cdot(-e_1,\epsilon,X_\cdot)
\]
 to denote the regeneration times in the $-e_1$ direction.
Assume that there are two opposite nonzero limiting velocities in directions $e_1$ and $-e_1$, i.e.,
\[
v_{e_1} \cdot v_{-e_1}>0.
\]
We let $\check{\mathbf P}(\cdot):=\mathbf{P}(\cdot|R_{-e_1}=\infty)$.

\begin{lemma}\label{l3}
Assume that there are two nonzero limiting velocities in direction $e_1$. We sample $(\epsilon,\tilde{X}_\cdot)$ according to $\hat{\mathbf P}$ and let
$\tilde{\tau}_\cdot=\tau_\cdot (e_1,\epsilon,\tilde{X}_\cdot)$ denote its regeneration times. 
For $n\ge 1$, we let 
\[
Y_\cdot^n=(Y_i^n)_{i=0}^{\tilde{\tau}_n}:=(\tilde{X}_{\tilde{\tau}_n-i})_{i=0}^{\tilde{\tau}_n}
\]
be the reversed path of $(\tilde{X}_i)_{i=0}^{\tilde{\tau}_n}$.
If $|z|$ is large enough, $d\ge 5$ and $n\ge 1$, then
\begin{equation}\label{e8}
E_{\hat{\mathbf P}}\check{\mathbf P}^{\tilde{X}_{\tilde{\tau}_n}+z}
\big(X_\cdot\in A(Y_\cdot^n, z)\big)>C>0.
\end{equation}
\end{lemma}

\pf
Let
\[m_z:=\lfloor|z|^{1/2}\rfloor.\]
Then
\begin{align}
&E_{\hat{\mathbf P}}\check{\mathbf P}^{\tilde{X}_{\tilde{\tau}_n}+z}
\big(X_\cdot\notin A(Y_\cdot^n, z)\big)\nonumber\\
&\le 
E_{\hat{\mathbf P}}\check{\mathbf P}^{\tilde{X}_{\tilde{\tau}_n}+z}(\tau'_{m_z}\ge |z|/3)+
\hat{\mathbf{P}}(\tilde{\tau}_n-\tilde{\tau}_{n-m_z}\ge |z|/3)\label{*12}\\
&\quad +E_{\hat{\mathbf P}}\check{\mathbf P}^{\tilde{X}_{\tilde{\tau}_n}+z}
(d(X_i,Y_\cdot^n)\le i^\delta\text{ for some } i>\tau'_{m_z})\label{*13}\\
&\quad +E_{\hat{\mathbf P}}\check{\mathbf P}^{\tilde{X}_{\tilde{\tau}_n}+z}
(d(\tilde{X}_{\tilde{\tau}_n-j},X_\cdot)\le j^\delta\text{ for some }j>\tilde{\tau}_n-\tilde{\tau}_{n-m_z}).\label{*14}
\end{align}

We will first estimate (\ref{*12}). By the translation invariance of the environment measure, 
\[
\check{\mathbf P}^x(\tau'_{m_z}\ge |z|/3)=\check{\mathbf P}(\tau'_{m_z}\ge |z|/3)
\text{ for any }x\in\mathbb{Z}^d.
\]
Hence
\begin{equation}\label{2e4}
E_{\hat{\mathbf P}}\check{\mathbf P}^{\tilde{X}_{\tilde{\tau}_n}+z}(\tau'_{m_z}\ge |z|/3)
=
\check{\mathbf P}(\tau'_{m_z}\ge |z|/3) 
\le 
\frac{3E_{\check{\mathbf P}}\tau'_{m_z}}{|z|}
\stackrel{(\ref{etau})}{\le}
C(E_{\check{\mathbf P}}\tau'_1)|z|^{-1/2}.
\end{equation}
Similarly,
\begin{equation}\label{2e5}
\hat{\mathbf P}(\tilde{\tau}_n-\tilde{\tau}_{n-m_z}\ge |z|/3)
\stackrel{(\ref{prop2})}{\le}
\exp{(e^{-cL})}\hat{\mathbf P}(\tau_{m_z}\ge |z|/3)\le C(E_{\hat{\mathbf P}}\tau_1) |z|^{-1/2}.
\end{equation}

To estimate (\ref{*13}) and (\ref{*14}), for $i\ge 1, n\ge j\ge 1$, we let
\begin{align*}
&Q'(i,x)=\check{\mathbf P}(x\text{ is visited in }[\tau'_{i-1},\tau'_i)),\\
&\tilde{Q}(j,x)=\hat{\mathbf P}(X_{\tau_n}+x \text{ is visited in}[\tau_{n-j},\tau_{n-j+1})).
\end{align*}
Note that by arguments that are similar to the proof of Theorem \ref{hke}, 
one can also obtain the heat kernel estimate
(\ref{ehke}) for $Q'(i,x)$ and $\tilde{Q}(j,x)$.
For $l>0$, let $B(o,l)=\{x\in\mathbb{Z}^d: d(o,x)\le l\}$.
Recall the definition of the $r$-boundary in Definition \ref{def1}.
By the translation invariance of the environment measure,
\[
\check{\mathbf P}^y(X_i=y+z)
=\check{\mathbf P}(X_i=z) \text{ for any }y,z\in\mathbb{Z}^d \text{and }i\in\mathbb{N}.
\]
Hence
\begin{align*}
&E_{\hat{\mathbf P}}\check{\mathbf P}^{\tilde{X}_{\tilde{\tau}_n}+z}
(d(X_i,\tilde{X}_\cdot)\le i^\delta\text{ for some } i>\tau'_{m_z})\\
&\le 
\sum_{i\ge m_z}\sum_{y\in\partial_1 B(o,i^\delta)}\sum_x
E_{\hat{\mathbf P}}\big[\check{\mathbf P}^{\tilde{X}_{\tilde{\tau}_n}+z}
(\tilde{X}_{\tilde{\tau}_n}+z+x\text{ is visited in }[\tau'_i,\tau'_{i+1}))\\
&\qquad\qquad\qquad\qquad\qquad\times 1_{\tilde{X}_{\tilde{\tau}_n}+z+x+y\in Y_\cdot^n}\big]\\
&=
\sum_{i\ge m_z}\sum_{y\in\partial_1 B(o,i^\delta)}\sum_x
\check{\mathbf P}(x\text{ is visited in }[\tau'_i,\tau'_{i+1}))
\hat{\mathbf P}(\tilde{X}_{\tilde{\tau}_n}+z+x+y\in Y_\cdot^n)\\
&=\sum_{i\ge m_z}\sum_{y\in\partial_1 B(o,i^\delta)}\sum_{j\le n}\sum_x
Q'(i,x)\tilde{Q}(j,x+z+y).
\end{align*}
By the heat kernel estimates and H\"{o}lder's inequality, 
\begin{align*}
\sum_{j\le n}\sum_x Q'(i,x)\tilde{Q}(j,x+z+y)
&\le 
\sqrt{\sum_x Q'(i,x)^2}\sum_{j\le n}\sqrt{\sum_x \tilde{Q}(j,x+y)^2}\\
&\le
 C(E_{\check{\mathbf P}}\tau'_1)i^{-d/4}
 \sum_{j\le n}(E_{\hat{\mathbf P}}\tau_1)j^{-d/4}\\
&\stackrel{d\ge 5}{\le} 
Ci^{-d/4}E_{\check{\mathbf P}}\tau'_1 E_{\hat{\mathbf P}}\tau_1.
\end{align*}
Thus
\begin{align}\label{2e6}
&E_{\hat{\mathbf P}}\check{\mathbf P}^{\tilde{X}_{\tilde{\tau}_n}+z}
(d(X_i,\tilde{X}_\cdot)\le i^\delta\text{ for some } i>\tau'_{m_z})\nonumber\\
&\le 
C\sum_{i\ge m_z}\sum_{y\in\partial_1 B(o,i^\delta)}i^{-d/4}
E_{\check{\mathbf P}}\tau'_1 E_{\hat{\mathbf P}}\tau_1\nonumber\\
&\le 
C\sum_{i\ge m_z}i^{(d-1)\delta}i^{-d/4}E_{\check{\mathbf P}}\tau'_1 E_{\hat{\mathbf P}}\tau_1
\le C 
|z|^{-(d-4)/8}
E_{\check{\mathbf P}}\tau'_1 E_{\hat{\mathbf P}}\tau_1,
\end{align}
where we used $d\ge 5$ and $\delta=\frac{d-4}{8(d-1)}$ in the last inequality.
Similarly, we have
\begin{align}\label{2e7}
&E_{\hat{\mathbf P}}\check{\mathbf P}^{\tilde{X}_{\tilde{\tau}_n}+z}
(d(\tilde{X}_{\tilde{\tau}_n-j},X_\cdot)
\le j^\delta\text{ for some }j>\tilde{\tau}_n-\tilde{\tau}_{n-m_z})\nonumber\\
&\le 
C |z|^{-(d-4)/8}E_{\check{\mathbf P}}\tau'_1 E_{\hat{\mathbf P}}\tau_1.
\end{align}

Combining (\ref{2e4}), (\ref{2e5}), (\ref{2e6}) and (\ref{2e7}), we conclude that
\begin{equation*}
E_{\hat{\mathbf P}}\check{\mathbf P}^{\tilde{X}_{\tilde{\tau}_n}+z}
\big(X_\cdot\in A(Y_\cdot^n, z)\big)>C>0,
\end{equation*}
if $|z|$ is large enough and $d\ge 5$.\qed\\

Let
\[
T^o=\inf\{i\ge 0: X_i\cdot e_1<0\}.
\]
For every fixed $\omega\in\Omega$ and $P_{\omega,\epsilon}^o$-almost every
$X_\cdot$,
\[
P_{\omega,\theta^n\epsilon}^{X_n}(T^o=\infty)1_{T^o>n}=P_{\omega,\epsilon}^o (T^o=\infty|X_1,\ldots,X_n),
\]
and so by Levy's martingale convergence theorem,
\[
\lim_{n\to\infty}P_{\omega,\theta^n\epsilon}^{X_n}(T^o=\infty)1_{T^o> n}= 1_{T^o=\infty}, 
\quad\text{$P_{\omega,\epsilon}^o$-almost surely}.
\]
Hence, for $(\omega, \epsilon,\tilde{X}_\cdot)$ sampled according to
$\hat{\mathbf P}$,
\[
\lim_{n\to\infty}
P_{\omega,\theta^{\tilde{\tau}_n}\epsilon}^{\tilde{X}_{\tilde{\tau}_n}}(T^o=\infty)=1, \quad\text{$\hat{\mathbf P}$-almost surely}.
\]
It then follows by the dominated convergence theorem that
\begin{equation}\label{3e7}
\lim_{n\to\infty}
E_{\hat{\mathbf P}}
P_{\omega,\theta^{\tilde{\tau}_n}\epsilon}^{\tilde{X}_{\tilde{\tau}_n}}(T^o<\infty)
=0.
\end{equation}

\begin{lemma}\label{l1}
For any $z\in\mathbb{Z}^d$,
\begin{equation}\label{e2}
\lim_{n\to\infty}
E_{\hat{\mathbf P}}
P_{\omega,\theta^{\tilde{\tau}_n}\epsilon}^{\tilde{X}_{\tilde{\tau}_n}+z}(T^o<\infty)
=0.
\end{equation}
\end{lemma}

\pf
For $n>|z|$, obviously
\[(\tilde{X}_{\tilde{\tau}_n}+z)\cdot e_1>0.\]
This together with ellipticity yields
\[
P_{\omega,\theta^{\tilde{\tau}_n}\epsilon}^{\tilde{X}_{\tilde{\tau}_n}}(T^o<\infty)
\ge{(\frac{\kappa}{2})^{|z|}}
P_{\omega,\theta^{\tilde{\tau}_n+|z|}\epsilon}^{\tilde{X}_{\tilde{\tau}_n}+z}(T^o<\infty).
\]
Hence using \eqref{3e7},
\[
\lim_{n\to\infty}
E_{\hat{\mathbf P}}
P_{\omega,\theta^{\tilde{\tau}_n+|z|}\epsilon}^{\tilde{X}_{\tilde{\tau}_n}+z}(T^o<\infty)
=0.
\]
On the other hand, noting that $\{R>\tau_1\}=\{R=\infty\}$,
\begin{align*}
&E_{\hat{\mathbf P}}
P_{\omega,\theta^{\tilde{\tau}_n+|z|}\epsilon}^{\tilde{X}_{\tilde{\tau}_n}+z}(T^o<\infty)\\
&= \sum_{m,x}
E_{P\otimes Q}[P_{\omega,\theta^{m+|z|}\epsilon}^{x+z}(T^o<\infty)
P_{\omega, \epsilon}^o (R>\tau_1,\tau_n=m,X_m=x)]/\mathbf{P}(R=\infty)\\
&=\sum_{m,x}
E_{P\otimes Q}[P_{\omega,\theta^{m}\epsilon}^{x+z}(T^o<\infty)
P_{\omega, \epsilon}^o (R>\tau_1,\tau_n=m,X_m=x)]/\mathbf{P}(R=\infty)\\
&=
E_{\hat{\mathbf P}}P_{\omega,\theta^{\tilde{\tau}_n}\epsilon}^{\tilde{X}_{\tilde{\tau}_n}+z}(T^o<\infty),
\end{align*}
where we used the independence (under $Q$) of $P_{\omega,
\theta^m\epsilon}^{x+z}(T^o<\infty)$ and
$P_{\omega,\epsilon}^o(R>\tau_1,\tau_n=m, X_m=x)$ in the second to last
equality. The conclusion follows.
\qed\\

\noindent\textit{Proof of the uniqueness of the non-zero velocity when $d\ge 5$, as stated in Theorem \ref{thm2}:} 
If the two antipodal velocities are both non-zero, we assume that 
\[v_{e_1}\cdot v_{-e_1}>0.\]

Sample $(\omega,\epsilon_\cdot,\tilde{X}_\cdot)$ according
to $\hat{\mathbf P}$.
Henceforth, we take $z=z_0$ such that (\ref{e8}) holds
and 
\[z_0\cdot e_1<-L.\] 


We will prove Theorem \ref{thm2} by showing that
\begin{equation}\label{contradiction}
E_{\hat{\mathbf P}} 
P_{\omega,\theta^{\tilde{\tau}_n}\epsilon}^{\tilde{X}_{\tilde{\tau}_n}+z_0}(T^o<\infty)>C
\end{equation}
for all $n>|z_0|$, which contradicts with (\ref{e2}).

First, let $\mathcal{G}$ denote the set of finite paths $y_\cdot=(y_i)_{i=0}^M$ that satisfy $y_M=0, M<\infty$.
Then
\begin{align}\label{e11}
& E_{\hat{\mathbf P}} 
P_{\omega,\theta^{\tilde{\tau}_n}\epsilon}^{\tilde{X}_{\tilde{\tau}_n}+z_0}(T^o<\infty)\\
&\ge 
E_{\hat{\mathbf P}} 
P_{\omega,\theta^{\tilde{\tau}_n}\epsilon}^{\tilde{X}_{\tilde{\tau}_n}+z_0}
\big((X_i)_{i=0}^{T^o}\in A(Y_\cdot^n, z_0),T^o<\infty\big)\nonumber\\
&=
\sum_{y_\cdot=(y_i)_{i=0}^M\in \mathcal{G}}
E_{\hat{\mathbf P}} 
[P_{\omega,\theta^M\epsilon}^{y_0+z_0}
\big((X_i)_{i=0}^{T^o}\in A(y_\cdot, z_0),T^o<\infty\big)
1_{Y_\cdot^n=y_\cdot}]\nonumber\\
&=
\frac{1}{\mathbf{P}(R=\infty)}
\sum_{y_\cdot\in \mathcal{G}}
\sum_{\substack{N<\infty\\(x_i)_{i=0}^N\in A(y_\cdot, z_0)}}
E_{P\otimes Q}
[P_{\omega,\theta^M\epsilon}^{y_0+z_0}
\big((X_i)_{i=0}^{T^o}=x_\cdot\big)
P_{\omega,\epsilon}(Y_\cdot^n=y_\cdot)].\nonumber
\end{align}
By the definition of the regeneration times, for any finite path $y_\cdot=(y_i)_{i=0}^M$,
there exists an event $G_{y_\cdot}$ such that $P_{\omega,\epsilon}(G_{y_\cdot})$ is
$\sigma(\epsilon_{i,y_i}, \omega_{y_j}:0\le i\le M, 0\le j\le M-L)$-measurable and
\[
\{Y_\cdot^n=y_\cdot\}=\{(\tilde{X}_i)_{i=0}^{\tilde{\tau}_n}=(y_{M-j})_{j=0}^M\}
=G_{y_\cdot}\cap\{R\circ\theta_M=\infty\}.
\]
Hence, for and any $y_\cdot=(y_i)_{i=0}^M\in \mathcal{G}$ and
$x_\cdot=(x_i)_{i=0}^N\in A(y_\cdot,z_0)$, $N<\infty$,
\begin{align}\label{2e8}
& E_{P\otimes Q}
[P_{\omega,\theta^M\epsilon}^{y_0+z_0}\big((X_i)_{i=0}^{T^o}=x_\cdot, R_{-e_1}>N\big)
P_{\omega,\epsilon}(Y_\cdot^n=y_\cdot)]\nonumber\\
&=E_P[
\bar{P}_\omega^{y_0+z_0}\big((X_i)_{i=0}^{T^o}=x_\cdot, R_{-e_1}>N\big)
\bar{P}_\omega(G_{y_\cdot})
\bar{P}_\omega^{y_0}(R=\infty)]\nonumber\\
&\stackrel{(\ref{prop1})}{\ge}
CE_P[
\bar{P}_\omega^{y_0+z_0}\big((X_i)_{i=0}^{T^o}=x_\cdot, R_{-e_1}>N\big)
\bar{P}_\omega(G_{y_\cdot})]
\bar{\mathbf P}(R=\infty).
\end{align}
where we used in the equality that $(\epsilon_{i,x})_{i\ge 0, x\in\mathbb{Z}^d}$
are iid and in the inequality the fact that
\[
\bar{P}_\omega^{y_0+z_0}\big((X_i)_{i=0}^{T^o}=x_\cdot, R_{-e_1}>N\big)
\bar{P}_\omega(G_{y_\cdot})
\]
is $\sigma(\omega_v: v\cdot e_1\le y_0\cdot e_1-L)$-measurable (note that $z_0\cdot e_1<-L$).
Further, by Lemma \ref{c2} and (\ref{e10}),
we have
\begin{align}\label{2e9}
&E_P[
\bar{P}_\omega^{y_0+z_0}\big((X_i)_{i=0}^{T^o}=x_\cdot, R_{-e_1}>N\big)
\bar{P}_\omega(G_{y_\cdot})]\nonumber\\
&\ge 
C\bar{\mathbf P}^{y_0+z_0}\big((X_i)_{i=0}^{T^o}=x_\cdot, R_{-e_1}>N\big)
\bar{\mathbf P}(G_{y_\cdot}).
\end{align}
Note that
\begin{align}\label{e12}
\bar{\mathbf P}(G_{y_\cdot})\bar{\mathbf P}(R=\infty)
&\stackrel{(\ref{prop1})}{\ge}
CE_P
[\bar{P}_\omega (G_{y_\cdot})\bar{P}_{\omega}^{y_0}(R=\infty)]
\nonumber\\
&=C\bar{\mathbf P}(Y_\cdot^n=y_\cdot)
\ge
C\hat{\mathbf P}(Y_\cdot^n=y_\cdot).
\end{align}
Therefore, by (\ref{e11}), (\ref{2e8}) and (\ref{2e9}),
\begin{align*}
& E_{\hat{\mathbf P}} 
P_{\omega,\theta^{\tilde{\tau}_n}\epsilon}^{\tilde{X}_{\tilde{\tau}_n}+z_0}(T^o<\infty)\\
&\ge
C\sum_{y_\cdot\in \mathcal{G}}
\sum_{\substack{N<\infty\\(x_i)_{i=0}^N\in A(y_\cdot, z_0)}}
\bar{\mathbf P}^{y_0+z_0}\big((X_i)_{i=0}^{T^o}=x_\cdot, R_{-e_1}>N\big)
\bar{\mathbf P}(G_{y_\cdot})\bar{\mathbf P}(R=\infty)\\
&\stackrel{(\ref{e12})}{\ge}
C\sum_{y_\cdot\in \mathcal{G}}
\sum_{\substack{N<\infty\\(x_i)_{i=0}^N\in A(y_\cdot, z_0)}}
\bar{\mathbf P}^{y_0+z_0}\big((X_i)_{i=0}^{T^o}=x_\cdot, R_{-e_1}>N\big)
\hat{\mathbf P}(Y_\cdot^n=y_\cdot)\\
&\ge 
CE_{\hat{\mathbf P}}
\check{\mathbf P}^{\tilde{X}_{\tilde{\tau}_n}+z_0}
(X_\cdot\in A(Y_\cdot^n, z_0))\stackrel{\text{Lemma }\ref{l3}}{>}C.
\end{align*}
Inequality (\ref{contradiction}) is proved.\qed

\textbf{Acknowlegments}.
I thank my advisor Ofer Zeitouni for his careful reading of earlier versions of this
paper and many valuable suggestions. Part of this work was done while I visited the 
Weizmann Institute of Science. I also thank the faculties and graduate students of
the math department for their hospitality. Finally, I am grateful to two anonymous referees, whose many insightful 
comments have led to a great improvement in the paper. 

This paper forms part of
my PhD thesis.


\begin{thebibliography}{99}
\bibitem{Be} N. Berger,
\emph{Limiting velocity of high-dimensional random walk in random environment},
Ann. Probab. 36 (2008), no. 2, 728-738.
\bibitem{CZ1}F. Comets, O.Zeitouni,
\emph{A law of large numbers for random walks in random mixing environments},
 Ann. Probab. 32 (2004), no. 1B, 880-914.
\bibitem{CZ2}F. Comets, O.Zeitouni,
\emph{Gaussian fluctuations for random walks in random mixing environments},
Probability in mathematics. Israel J. Math. 148 (2005), 87-113. 
\bibitem{DS}R. Dobrushin and S. Shlosman,
\emph{Completely analytical Gibbs fields},
Statistical physics and dynamical systems (K\"oszeg, 1984), 371-403,
Progr. Phys., 10, Birkh\"auser Boston, Boston, MA, 1985.
\bibitem{Go} L. Goergen,
\emph{Limit velocity and zero-one laws for diffusions in random environment},
Ann. Appl. Probab. 16 (2006), no. 3, 1086-1123. 
\bibitem{R-A1}F. Rassoul-Agha,
\emph{The point of view of the particle on the law of large numbers for random walks in a mixing random environment},
Ann. Probab. 31 (2003), no. 3, 1441-1463. 
\bibitem{R-A2}F. Rassoul-Agha,
\emph{Large deviations for random walks in a mixing random environment and other (non-Markov) random walks},
Comm. Pure Appl. Math. 57 (2004), no. 9, 1178-1196.
\bibitem{R-A3}F. Rassoul-Agha,
\emph{On the zero-one law and the law of large numbers for random walk in mixing random environment},
Electron. Comm. Probab. 10 (2005), 36-44.
\bibitem{Tho}H. Thorisson,
\emph{Coupling, stationarity, and regeneration},
Probability and its Applications (New York). Springer-Verlag, New York, 2000.
\bibitem{So}F. Solomon,
\emph{Random walks in a random environment},
Ann. Probability 3 (1975), 1-31.
\bibitem{SZ}A.S. Sznitman, M. Zerner,
\emph{A law of large numbers for random walks in random environment},
Ann. Probab. 27 (1999), no. 4, 1851-1869.
\bibitem{ZO} O. Zeitouni,
\emph{Random walks in random environment}, Lectures on probability theory and statistics, 189-312,
Lecture Notes in Math., 1837, Springer, Berlin, 2004. 
\bibitem{Ze} M. Zerner,
\emph{A non-ballistic law of large numbers for random walks in i.i.d. random environment},
Electron. Comm. Probab. 7 (2002), 191-197.
\bibitem{ZM} M. Zerner, F. Merkl,
\emph{A zero-one law for planar random walks in random environment},
Ann. Probab. 29 (2001), no. 4, 1716-1732.
\end{thebibliography}
\end{document}